\documentclass[review]{elsarticle}

\usepackage{lineno}
\usepackage{theorem}
\usepackage{amssymb,amsmath}
\usepackage{bbm}
\usepackage{bm}

\usepackage{mathptmx}       
\usepackage{helvet}         
\usepackage{courier}        
\usepackage{type1cm}        
\usepackage{makeidx}         
\usepackage{graphicx}        
\usepackage{multicol}        
\usepackage[bottom]{footmisc}
\usepackage{subfigure}
\usepackage{booktabs}

\usepackage{epsfig} 
\usepackage{epstopdf}
\usepackage{color}
\usepackage{colordvi}

\newcommand{\duality}[4]         { \big \langle {#1},{#2} \big \rangle_{{#3}\times{#4}} }       
\newcommand{\dualp}[2]         { \big \langle {#1},{#2} \big \rangle }       
\newcommand{\norm}[2]         { \| {#1} \|_{#2} }                      
\newcommand{\bfm}[1]             { \mathbf{#1}     }             %
\newcommand{\rr}             { \bfm{r}     }             
\newcommand{\boldphi}             { \boldsymbol{\varphi}     }             
\newcommand{\Nabla}       { \boldsymbol{\nabla} }   
\newcommand{\SLTO}            { L^2(\Omega) }                       
\newcommand{\SLTOT}            { L^2(\Omega_T) }                       
\newcommand{\SLTOvec}            { [L^2(\Omega)]^2 }                   
\newcommand{\SHOO}            { H^1(\Omega) }                         
\newcommand{\SHOP}            { H^1(\Ph) }                              
\newcommand{\SHOPT}            { H^1(\textcolor{black}{\PhT}) }   
\newcommand{\SHdivP}            { H(\text{div},\Ph) }            
\newcommand{\VVK}            { {V(K_m)} }              
\newcommand{\VOM}            { {V(\Omega)} }              
\newcommand{\VV}            { {V(\Ph)} }        
\newcommand{\VVT}            { {V(\textcolor{black}{\PhT})} }        
\newcommand{\UUU}            { U(\Omega) }           
\newcommand{\UUUA}            { U^A(\Omega) }           
\newcommand{\UUUhat}            { \hat{U}(\Gamma_h) }           
\newcommand{\UUUT}            { U(\Omega_T) }           
\newcommand{\UUUh}            { U^h(\Omega) }           
\newcommand{\UUUhT}            { U^h(\Omega_T) }           
\newcommand{\VVh}            { V^*(\Ph) }          
\newcommand{\VVhT}            { V^*(\PhT) }          
\newcommand{\SHOK}            { H^1(K_m) }                        
\newcommand{\SHdivO}            { H(\text{div},\Omega) }             
\newcommand{\SHdivK}            { H(\text{div},K_m) }            
\newcommand{\SHHdK}            { H^{1/2}(\dKm) }         
\newcommand{\SHmHdK}            { H^{-1/2}(\dKm)  }       
\newcommand{\SHOOT}            { H^1(\Omega_T) }                   
\newcommand{\SLOinf}            { L^\infty(\Omega)}          

\newcommand{\xx}             { \bfm{x}     }             
\newcommand{\vn}             { \bfm{n}     }             

\newcommand{\Ph}            { \mathcal{P}_h }
\newcommand{\PhT}            { \textcolor{black}{\mathcal{P}_T^h} }
\newcommand{\Kep}           { K_m \in \Ph}
\newcommand{\KepT}           { K_m \in \PhT}
\newcommand{\dKm}           { \partial K_m }

\newcommand{\dx}              { \; {\rm d} \bfm{x}   }                
\newcommand{\dss}             { \, {\rm d} s   }                          
\newcommand{\summa}[2]        { \overset{#2}{\underset{#1}{\sum}} } 
\newcommand{\supp}[1]         { \underset{#1}{\sup} \, }        

\newcommand{\sss}              { \bfm{s} }                      
\newcommand{\cc}              { \bfm{c} }                      
\newcommand{\dd}              { \bfm{d} }                      
\newcommand{\nn}              { \bfm{n} }                      
\newcommand{\pp}              { \bfm{p} }       
\newcommand{\ttt}             { \bfm{t}     }             
\newcommand{\eps}             { \boldsymbol{\varepsilon} }     
\newcommand{\isdef}           { \overset{\text{def}}{=} } 
\newcommand{\ds}              { \displaystyle }   

\theoremstyle{plain}
\theoremheaderfont{\bfseries \upshape}

\newtheorem{lem}{Lemma}[section]
\newtheorem{rem}{Remark}[section]

\newtheorem{prp}{Proposition}[section]

\usepackage[margin=1in]{geometry}

\begin{document}

\begin{frontmatter}
 \title{A Stable FE Method For the Space-Time Solution of the Cahn-Hilliard Equation}

\author[1]{Eirik Valseth\corref{cor1}}
\ead{Eirik@utexas.edu}

\author[2]{Albert Romkes}
\ead{Albert.Romkes@sdsmt.edu}

\author[2]{Austin R. Kaul}
\ead{Austin.Kaul@mines.sdsmt.edu}


 \cortext[cor1]{Corresponding author}

 \address[1]{Oden Institute for Computational Engineering and Sciences, The University of Texas at Austin, Austin, TX 78712, USA}
 \address[2]{Department of Mechanical Engineering, South Dakota School of Mines \& Technology, Rapid City, SD 57701, USA}


\biboptions{sort&compress}


%
%
%
\begin{abstract}
In its application to the modeling of a mineral separation process, we propose the numerical analysis of the Cahn-Hilliard equation by employing space-time 
discretizations of  
the automatic variationally stable finite element (AVS-FE) method.
 The AVS-FE method is a Petrov-Galerkin method which employs the concept of optimal discontinuous test functions of the 
discontinuous Petrov-Galerkin (DPG) method by Demkowicz and Gopalakrishnan. The trial space, however, consists of globally continuous Hilbert spaces 
such as $\SHOO$ and $\SHdivO$. Hence, the AVS-FE approximations employ classical $C^0$ or Raviart-Thomas FE basis functions.
The optimal test functions guarantee the numerical stability of the AVS-FE method and lead to 
discrete systems that are symmetric and positive definite. Hence, the AVS-FE method can 
solve the Cahn-Hilliard equation in both space and time without a restrictive CFL 
condition to dictate the space-time element size. We present \textcolor{black}{multiple} numerical verifications of both \textcolor{black}{stationary and transient problems}. The verifications show optimal rates of 
convergence in $\SLTO$ and $\SHOO$ norms. 
Results for mesh adaptive refinements \textcolor{black}{in both space and time} using a 
built-in error estimator of the AVS-FE method are also presented.  
\end{abstract}

\end{frontmatter}

%
\section{Introduction}
\label{sec:introduction}

The refinement and concentration of minerals from mineral ores is a process that 
typically requires the use of water such as flowing film and froth flotation concentrators.
Processing facilities in the United States consume large amounts of water, in some cases up to $60,000 m^3$ each day~\cite{mineralconf}.
 Therefore, sustainable approaches to mineral 
concentration that significantly reduce or remove the need for water are needed to 
aid in conservation efforts. Furthermore, the location of several copper mines in the United States are in arid regions of the Southwest, thereby further increasing the 
importance of conservation efforts. 
It has been proposed by researchers at South Dakota School of Mines \& Technology 
(SDSM\&T) to exploit the adhesion forces between mineral particles and specifically 
tailored substrates to develop new mineral separation techniques using as little 
water as possible.  
Thus, a new type of mineral separator must be developed and designed. To aid in the 
design process, it is necessary to predict the mineral separation process which  requires simulation of the accumulation of mineral particles on chemically treated substrates. This accumulation is to be modeled by the Cahn-Hilliard equation.


The mathematical analysis and well posedness results for the Cahn-Hilliard equation have been established in, e.g.,~\cite{elliott1986cahn} by Elliott and Zheng, thereby setting the stage for the application of a FE method in its approximation.  However, there 
are two challenges:
$i)$ the nonlinearity of the Cahn-Hilliard equation and $ii)$ the transient nature of this problem leading to a loss of the numerical stability for the \textcolor{black}{Galerkin} FE method.  
The second challenge is typically critical, as FE methods for nonlinear problems have been established
successfully in, e.g.,~\cite{oden2006finite}.
\textcolor{black}{To achieve stability in the classical \textcolor{black}{Galerkin} FE method, the FE mesh partition, i.e., element size, must be fine enough to establish numerically stable FE approximations thereby leaving the Galerkin FE method unsuitable for a space-time approximation of the  Cahn-Hilliard equation.} 
\textcolor{black}{When solving transient partial differential equations (PDE)s in a FE framework, a method of lines approach is typically taken, i.e.,  FE methods are employed in the spatial domain whereas the temporal domain is discretized by a 
finite difference scheme.} The numerical stability of finite difference schemes is then established through the Courant-Friedrichs-Lewy (CFL) condition~\cite{courant1928partiellen}. 
This approach has been employed by several authors to establish approximations of 
the Cahn-Hilliard equation using several flavors of FE methods including least squares FE method, discontinuous Galerkin methods, and isogeometric analysis (see~\cite{chave2016hybrid,clavijo2019reactive,wells2006discontinuous,
barrett2001fully,gomez2008isogeometric,fernandino2011least,dean1996approximate}). 

An alternative to temporal discretizations using  difference methods and a CFL condition 
are conditionally stable FE methods which have been successfully applied to transient problems, see, 
e.g.,~\cite{Hughes1996,hughes1988space,aziz1989continuous}. While these space-time methods have been
 successful in the FE approximation of transient phenomena, their conditional stability requires arduous 
a priori analyses to properly determine their stabilization parameters.
\textcolor{black}{Guaranteed} stable FE methods are also applicable for transient PDEs such as the DPG method~\cite{ellis2014space,ellis2016robust,roberts2015discontinuous} or least squares FE methods~\cite{bochevLeastSquares}.
Fernandino and Dorao~\cite{fernandino2011least} applied the least squares FE method successfully to a Cahn-Hilliard problem in both space and time using basis functions that are of higher 
order continuity than classical \textcolor{black}{Galerkin} FE methods.
The computational cost of these space-time FE methods are typically higher than the method of lines approach but the 
FE formulations \textcolor{black}{have the advantage that they} can employ a wide range of tools such as \emph{a priori} and \emph{a posteriori} error 
estimation and $hp-$adaptive refinement strategies in both space and time. Thus, the (potential) additional computational cost can be justified.

The AVS-FE method, introduced by Calo, Romkes, and Valseth in\textcolor{black}{~\cite{CaloRomkesValseth2018}}, is a stable FE 
method, i.e., the AVS-FE approximations are guaranteed to remain stable for 
any PDE as long as the kernel of the underlying differential operator is trivial \textcolor{black}{and the optimal test functions are resolved with sufficient accuracy}. This method is a Petrov-Galerkin method in which the trial space consist of globally continuous FE bases and the test space of piecewise discontinuous functions. Hence, it is a hybrid between the DPG method of Demkowicz and Gopalakrishnan~\cite{Demkowicz4} and classical \textcolor{black}{Galerkin} FE methods. In addition to its discrete stability, other features of the AVS-FE method are 
highly accurate flux approximations, due to its first-order system setting, and its ability to compute optimal 
discontinuous test functions on the fly, element-by-element. \textcolor{black}{Other related methods are the first-order system least squares FE method~\cite{bochevLeastSquares} and the method of Calo \emph{et al.}~\cite{calo2019adaptive} in which a discretely stable discontinuous Galerkin formulation is used in a minimum residual setting.}

In this paper, we develop space-time AVS-FE approximations of the Cahn-Hilliard equation. 
A space-time AVS-FE method is chosen to exploit its stability property, 
its convergence properties, and the built-in error indicators allowing us to employ mesh adaptive 
refinement strategies.
The Cahn-Hilliard equation being nonlinear requires special treatment and we 
take the approach of Carstensen \emph{et al.} in~\cite{carstensen2018nonlinear}.
To start, we introduce the mineral processing application and the corresponding  model Cahn-Hilliard boundary value problem (BVP), in addition to notations and conventions in Section~\ref{sec:model_and_conv}.
Next, we review the AVS-FE methodology in Section~\ref{sec:AVS-review} for linear problems and introduce the concepts
of~\cite{carstensen2018nonlinear} to be employed to perform nonlinear iterations.
 In Section~\ref{sec:avs-fe}, we derive the equivalent AVS-FE weak 
formulation for the Cahn-Hilliard BVP.
In Section~\ref{sec:verifications}, we  perform multiple  numerical verifications.  First, in Section~\ref{sec:1d_prob},  we 
present verifications \textcolor{black}{for problems with manufactured exact solutions to assess convergence properties under both uniform and adaptive mesh refinements.} 
Then, in Section~\ref{sec:2d_prob}, we present a phase separation problem for the Cahn-Hilliard problem 
from literature that is spatially two dimensional. \textcolor{black}{The last numerical verification we 
consider is a heuristic model for a mineral separation process to verify the applicability of 
the Cahn-Hilliard equation to mineral separation. }
Finally, we conclude with remarks on the results and future 
works in Section~\ref{sec:conclusions}.

\section{Model Problem and Review of the AVS-FE Method}  
\label{sec:model_and_conv}


\subsection{Model Problem}
\label{sec:model_problem}

The Cahn-Hilliard equation~\cite{cahn1958free} was introduced to model the evolution of the phase transition of components in 
a binary alloy  from a mixed to a separated state.
 The equation is a fourth order 
nonlinear PDE and can be found in several forms in literature. Here, we consider the following form:
\begin{equation}  \label{eq:Cahn-Hilliard_eq}
 \ds - \frac{\partial u}{\partial t} \, +\, D \, \Delta \left[ u^3 - u   -\lambda \, \Delta \, u \right]   = 0,
\end{equation}
where \textcolor{black}{$u=u(\xx,t)$ denotes the concentration of a constituent undergoing a phase transition}, $\Delta$ is the spatial Laplacian, $D \in \SLOinf$ is the diffusion coefficient, and
$\lambda \in \SLOinf$ is the square of the width of the transition region in the separation
process. Note that $D$ is of unit $\frac{m^2}{s}$ and $\sqrt{\lambda}$ is of unit
$m^2$. 

\begin{figure}[h]
{\centering
\hspace{0.9in} \includegraphics[width=0.65\textwidth]{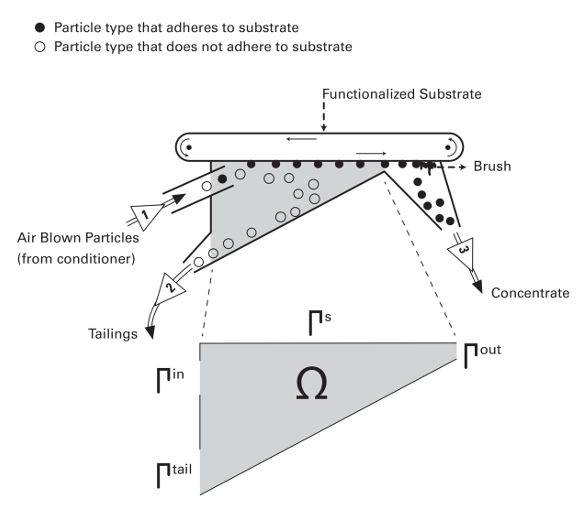}}
  \caption{\label{fig:mineral_separator} Conceptual sketch of mineral separator (Courtesy of Brian Hill of University Relations at SDSM\&T).}
\end{figure}

In collaboration with a team of metallurgists at the SDSM\&T, we propose to use the Cahn-Hilliard equation to model a particular mineral separation process. 
The proposed mineral separation process will exploit the forces of adhesion between mineral particles and substrates, both of
which are potentially treated, i.e., functionalized to ensure maximum adhesion of desirable minerals. 
In Figure~\ref{fig:mineral_separator}, a conceptual sketch of the process
is shown. Ore enters the separator in a stream of air, where the 
desirable minerals adhere to the substrate and the remainder of the ore ends up in what is referred to as tailings. 
The motivation for the development of this type of separation process stems from the fact 
that it would greatly reduce the use of water compared to currently used techniques. These generally 
consist of flotation processes which require large amounts of water to establish proper mineral 
separation from the ore.  Now, in the design of the newly proposed separation method, it is 
necessary to predict the separation of the mineral concentration as it accumulates onto the substrate.
To do so, we use the Cahn-Hilliard equation. 

In the particular application of the Cahn-Hilliard equation to the mineral separation process, the concentration function
$u = u(\xx,t)$ represents the mineral concentration at a material point $\xx$ and time $t$. 
Thus, in the spirit of the original application of this equation, a value 
$u(\xx,t) = 1$ represents the scenario in
which the desired mineral has separated from the flow of minerals and has adhered to the substrate. Conversely, a value 
$u(\xx,t) = 0$ denotes the scenario in which the desired mineral particle is still fully dispersed in the flow of minerals.
Values between 0 and 1 identify areas in which the mineral
particles are in the process of separating.

For the mathematical model of the separation process, we consider only the portion of the mineral separator 
in which the mineral separation occurs.
The boundary, which encloses the separator 
is assumed to consist of several disjoint portions corresponding to the substrate (onto 
which the mineral accumulates), the inflow boundary (where the minerals enter the 
separator), the tailing boundary (where minerals that fail to adhere to the substrate exit the separator), and finally the outflow boundary (where the accumulated minerals exit the separator). Additionally, there may be additional portions of the boundary that serve
to encompass the separation process. In Figure~\ref{fig:mineral_separator}, an example mineral separator computational domain $\Omega$ is shown.

To establish approximations of the Cahn-Hilliard equation, proper boundary and initial conditions are needed to guarantee physically meaningful results.
Let $\Omega \subset \mathbb{R}^2$ be an open bounded domain
(see Figure~\ref{fig:mineral_separator}) with Lipschitz boundary $\partial \Omega$
and outward unit normal vector $\vn$. Also, let $t = 0$ be the initial time and $t = T$ the final time.
The boundary $\partial \Omega$ consists of  open
subsections:
\begin{itemize}
\item $\Gamma_{s}$ - the portion of $\partial \Omega$ that coincides with the functionalized substrate upon which mineral accumulation occurs. 
\item $\Gamma_{in}$ - the inflow boundary.
\item $\Gamma_{out}$ - the portion of the boundary $\partial \Omega$ through which the accumulated separated mineral exits the separator.
\item $\Gamma_{tail}$ - the outflow boundary, containing the mineral particles that failed to adhere.
\item $\Gamma_{0}$ - the remaining portion of the boundary $\partial \Omega$ which contains the entire separation process.
\end{itemize}
The intersection of these portions is empty and $\partial \Omega = \overline{\Gamma_{s}\cup\Gamma_{in}\cup\Gamma_{out}\cup\Gamma_{tail}\cup\Gamma_{0}}$.
With the boundary identified by these subsets, the 
 proper boundary conditions for the case of mineral separation are as follows:
\begin{itemize}
\item $u = 1, \, \xx \in \Gamma_{s},\, t \in (0,T)$, i.e., the mineral particles have adhered to the substrate and separated from the flow field .
\item $u = 0, \, \xx \in \partial \Omega \setminus \Gamma_{s},\, t \in (0,T)$, i.e., no mineral particles adhere to the remainder of the separator.
\item $u^3 - u -\lambda \, \Delta \, u = 0, \, \xx \in \partial \Omega, \, t \in (0,T)$, if there is no separation ongoing, the "chemical potential" $u^3 - u -\lambda \, \Delta \, u$ must vanish .
\end{itemize}
Finally, the initial conditions for the mineral separation application are:
\begin{itemize}
\item $u = u_{in}, \, \xx \in \Gamma_{in},\, t =0$, i.e., we assume the separator has no material in it at the onset of the 
separation process. Hence, the initial conditions are zero on all boundaries.
\item $u = 0, \, \xx \in \partial \Omega \setminus \Gamma_{s},\, t = 0$.
\item $u = 0, \, \xx \in  \Omega,\; t = 0$.
\end{itemize}
where $u_{in}$ is the concentration of the desired mineral as it enters the separator.

With boundary and initial conditions defined along with the PDE~\eqref{eq:Cahn-Hilliard_eq}, we consider the following Cahn-Hilliard initial boundary value problem  (IBVP):
\begin{equation}  \label{eq:Cahn-Hilliard_BVP}
\boxed{
\begin{array}{l}
\text{Find }  u  \text{ such that:}    
\\[0.05in] 
\qquad 
\begin{array}{rcl}
 \ds - \frac{\partial u}{\partial t} \, +\, D \, \Delta \left[ u^3 - u  
 -\lambda \, \Delta \, u \right]  & = & 0, \quad \text{ in } \, \Omega\times(0,T), 
 \\[0.05in]
 \qquad \ds u &  = & u_0, \quad \text{ on } \,  \partial \Omega , 
 \\
  \qquad \ds u &  = & u_{initial}, \quad \text{ in } \, \textcolor{black}{ \Omega} , 
 \\
 \qquad u^3 - u -\lambda \, \Delta \, u & = & g, \quad \text{ on } \,  \partial \Omega, \\
 \end{array}
 \end{array}
}
\end{equation}
where the values of $u_0$, $u_{initial}$, and $g$ are as given in the  preceding lists. For generality, 
we keep these arbitrary and employ specific choices in the numerical verifications of Section~\ref{sec:verifications}.

While the goal of the current research project is to use the Cahn-Hillard equation as 
a model problem for mineral separation and the design of a mineral separator, we limit our presentation to the numerical 
approximation of the Cahn-Hilliard initial boundary value problem (IBVP) on general domains $\Omega$ as this is a key stepping stone towards the research goal and the design of the separator is still work in progress.

\subsection{Review of the AVS-FE method}
\label{sec:AVS-review}
The AVS-FE method has been introduced introduced by Calo, Romkes and Valseth in~\cite{CaloRomkesValseth2018}. 
 Attractive features of the AVS-FE method \textcolor{black}{to be exploited for the Cahn-Hilliard problem} are its \textcolor{black}{discrete} numerical 
stability property \textcolor{black}{and ease of adaptive mesh refinements due to its built-in error estimator and indicators}.
 The AVS-FE method 
is a Petrov-Galerkin method in which the 
trial space consists of continuous FE basis functions, and the 
test space consists of piecewise discontinuous functions. The discontinuous test space 
is spanned by so called 'optimal' test functions that are computed on-the-fly by invoking the DPG philosophy~\cite{Demkowicz4, Demkowicz2, Demkowicz3, Demkowicz5, Demkowicz6}. 

To introduce the AVS-FE method here, we consider a domain $\Omega$ partitioned into
elements:
\begin{equation}
\notag
\label{eq:domain}
  \Omega = \text{int} ( \bigcup_{\Kep} \overline{K_m} ),
\end{equation}
and an abstract AVS-FE weak form in which the underlying differential operator is linear:
\begin{equation} \label{eq:abstract_weak_form_AVS}
\boxed{
\begin{array}{ll}
\text{Find } \mathbbm{u} \in \UUU & \hspace{-0.05in} \text{ such that:}
\\[0.05in]
 &  \quad B(\mathbbm{u},\mathbbm{v}) = F(\mathbbm{v}), \quad \forall \mathbbm{v} \in \VV, 
 \end{array}}
\end{equation}
where $\mathbbm{u}$ and $\mathbbm{v}$ are the vector valued trial and test functions, respectively, $\UUU$ is the trial space, $\VV$ the broken test
space, $B:\UUU\times\VV\longrightarrow \mathbb{R}$ is the bilinear form,
  $F:\VV\longrightarrow \mathbb{R}$  the linear 'load' functional, and $\Ph$ denotes the partition of $\Omega$ into finite  elements (see~\eqref{eq:domain}).
In the AVS-FE method, $\UUU$ is a globally continuous Hilbert space as used in \textcolor{black}{mixed and Galerkin} FE methods.
However, $\VV$ is a broken space consisting of functions that are globally in $\SLTO$  
and locally may be of higher order (e.g., $\SHOK$). The kernel of $B(\cdot,\cdot)$ is assumed to be trivial to guarantee
the uniqueness of solutions (as in any other FE method). 

With the assumption on the kernel of the bilinear form in place,
the AVS-FE method introduces the following \emph{energy norm} $\norm{\cdot}{B}: \UUU\longrightarrow [0,\infty)$:
\begin{equation}
\label{eq:abstract_energy_norm}
\norm{\mathbbm{u}}{B} \isdef \supp{\mathbbm{v}\in \VV\setminus \{\mathbf{0}\}} 
     \frac{|B(\mathbbm{u},\mathbbm{v})|}{\norm{\mathbbm{v}}{\VV}}.
\end{equation}
%
%
%
%
%
%
%
%
The well posedness of the AVS-FE weak formulation is then established by the following lemma:
\begin{lem}
\label{lem:well_posed_cont}
Let the source and Neumann data \textcolor{black}{(if present)} be sufficiently regular. Then, the weak formulation~\eqref{eq:abstract_weak_form_AVS} is well posed.
\end{lem} 
\emph{Proof}: The proof follows from the Generalized Lax-Milgram Theorem, as $B(\cdot,\cdot)$ 
satisfies the Inf-sup condition as well as the continuity condition in terms of the energy
norm~\eqref{eq:abstract_energy_norm}, (see~\cite{Demkowicz5,Demkowicz4} for details).
\newline \noindent ~\qed

By deriving the weak statement such that the trial space \textcolor{black}{consists of} global Hilbert spaces,  
the AVS-FE method seeks FE approximations $\mathbbm{u}^h$ of $\mathbbm{u}$ 
 of~\eqref{eq:abstract_weak_form_AVS} in which
the trial functions in the discretization are
FE basis functions that span the FE trial space $\UUUh$, e.g., $\SHOO$ or $\SHdivO$.
Hence, we represent the approximations of the components  $\mathbbm{u}^h (\xx)$ as linear 
combinations of trial \textcolor{black}{basis} functions $e^i(\xx)\in\UUUh$ and the corresponding degrees of freedom, $\mathbbm{u}^h_i$.
Conversely, to construct the test space $\VVh$ we compute piecewise discontinuous \emph{optimal} test functions that
guarantee stable discretizations. 
These optimal test functions are obtained by employing 
the DPG philosophy~\cite{Demkowicz4, Demkowicz2, Demkowicz3, Demkowicz5, Demkowicz6} in which \emph{global} optimal test functions are established through  
\emph{global} weak problems. However, even  though the optimal test functions are global functions,
they have compact support and, in the case of the AVS-FE method their support is identical to that 
of the trial functions. Additionally, the local restrictions are computed in a completely 
decoupled fashion, i.e., element-by-element, with high accuracy (see~\cite{CaloRomkesValseth2018} for 
details).
Thus, e.g., for  the local restriction of a trial function $\mathbbm{e}^i$ on an 
element $\Kep,$ i.e., a
 shape function, we solve the corresponding optimal test function 
$\hat{\mathbbm{e}}^i(\xx)$ from the following \emph{local} problem on $K_m$:
\begin{equation}
\label{eq:local_test_problems}
\begin{array}{rcll}
\ds \left(\, \mathbbm{r},\hat{\mathbbm{e}}^i \, \right)_\VVK & = & B_{|K_m}(\,\mathbbm{e}^i,\mathbbm{r} \, ),& \quad \forall \mathbbm{r}\in\VVK, 
\end{array}
\end{equation}
where $B_{|K_m}(\cdot,\cdot)$ denotes the restriction of  $B(\cdot,\cdot)$ to the element
$K_m$, $\VVK$ the local restriction of the test space to $K_m$, and 
 $\left(\, \cdot, \cdot \, \right)_\VVK:\; \VVK\times\VVK \longrightarrow \mathbb{R}$, is a 
 local inner product on $\VVK$. 
see~\cite{CaloRomkesValseth2018} for details.
\textcolor{black}{Numerical evidence suggests that the local Riesz representation 
problems~\eqref{eq:local_test_problems} can be solved at the same local order of approximation as 
the trial function in the RHS of~\eqref{eq:local_test_problems}. Since this space consists of
discontinuous polynomial functions, it is larger than the space of the continuous trial
functions, i.e., this is in line with the DPG method of testing with a larger space to attain
discrete stability.  }

Finally, we introduce the FE discretization of~\eqref{eq:abstract_weak_form_AVS}  governing the approximation $\mathbbm{u}^h\in\UUUh$ of $\mathbbm{u}$ :
\begin{equation} \label{eq:discrete_form}
\boxed{
\begin{array}{ll}
\text{Find} &  \mathbbm{u}^h \in \UUUh \; \text{ such that:}
\\[0.1in]
 &   B(\mathbbm{u}^h,\mathbbm{v}^h) = F(\mathbbm{v}^h), \quad \forall \mathbbm{v}^h\in \VVh, 
 \end{array}}
\end{equation}
where the finite dimensional subspace of test functions $\VVh\subset\VV$ is spanned by the  optimal test functions.

By using the DPG philosophy to construct $\VVh$, the
discrete problem~\eqref{eq:discrete_form} inherits the continuity and inf-sup constants
of the continuous problem \textcolor{black}{scaled by the continuity constant of a Fortin type operator~\cite{nagaraj2017construction} }.
Hence, the AVS-FE discretization is stable for any choice of element size $h_m$ and local degree of polynomial approximation $p_m$. A further consequence of the optimal test functions is that the
global stiffness matrix is symmetric and positive definite regardless of the character of the underlying differential operator.

\begin{rem}
Instead of computing the optimal test functions from~\eqref{eq:local_test_problems} on-the-fly to construct
the FE system of linear algebraic equations, one can 
consider another, equivalent, interpretation of the AVS-FE method. 
This interpretation is in the DPG
literature~\cite{demkowicz2014overview,wozniak2016fast,roberts2014dpg}  referred to as a mixed or saddle point problem and is a result of the fact that DPG and AVS-FE methods are constrained 
minimization techniques:.
\begin{equation} \label{eq:constrained_min_problem}
\boxed{
\begin{array}{rl}
\text{Find } \mathbbm{u}^h \in  \UUUh, \mathbbm{\hat{E}}^h \in V^h(\Ph)  & \hspace{-0.15in} \text{ such that:}
\\[0.05in]
   \quad \left(\, \mathbbm{\hat{E}}^h,\mathbbm{v}^h \, \right)_\VV - B(\mathbbm{u}^h,\mathbbm{v}^h) & =  - F(\mathbbm{v}^h), \quad \forall \mathbbm{v}^h \in V^h(\Ph),  \\
  \quad B(\mathbbm{p}^h,\mathbbm{\hat{E}}^h)& =  0, \quad \forall \, \mathbbm{p}^h \in \UUUh.
 \end{array}}
\end{equation}
The second equation of~\eqref{eq:constrained_min_problem} represents the 
constraint in which the Gateaux derivative of the bilinear form is acting on the approximate
"error representation" function $\mathbbm{\hat{E}}^h$. This function is 
a Riesz representer of the approximation error $\mathbbm{u}-\mathbbm{u}^h$ and leads to
 \textcolor{black}{an}
identity between the energy norm of the approximation error and the norm of the error representation function on $\VV$. Hence, the 
norm $\norm{\mathbbm{\hat{E}}}{\VV}$ is an a posteriori error estimate and its local restriction may be employed as an error indicator in mesh adaptive strategies. For details on these error indicators and the derivation of the mixed formulation, see~\cite{demkowicz2014overview} or 
\textcolor{black}{\cite{cohen2012adaptivity}}.
Note that since we, at this point, have assumed that the underlying differential 
operator is linear, the Gateaux derivative of the bilinear form is identical 
to itself.

This mixed form allows straightforward implementation in high level FE solvers such as Firedrake~\cite{rathgeber2017firedrake} and FEniCS~\cite{alnaes2015fenics}. The  cost of solving the resulting system of linear algebraic equations from~\eqref{eq:constrained_min_problem} is larger than the 'classical' AVS-FE method since now the optimal test functions are essentially computed by solving global problems. However, it has the clear advantage for $hp$-adaptive strategies, since upon 
solving~\eqref{eq:constrained_min_problem}, it immediately provides \emph{a posteriori} error 
estimators and error indicators that can drive the mesh adaptive process.  
\end{rem}

\section{AVS-FE Weak Formulation and Discretization of The Cahn-Hilliard Equation}  
\label{sec:avs-fe}
With the notations introduced in Section~\ref{sec:model_and_conv} and the review of the
AVS-FE method above, we proceed to derive the AVS-FE weak formulation for the 
Cahn-Hilliard IBVP. To this end, let us consider the following general form of the 
Cahn-Hilliard IBVP~\eqref{eq:Cahn-Hilliard_BVP}:
\begin{equation}  \label{eq:gen_Cahn-Hilliard_BVP}
\boxed{
\begin{array}{l}
\text{Find }  u  \text{ such that:}    
\\[0.05in] 
\qquad 
\begin{array}{rcl}
 \ds - \frac{\partial u}{\partial t} \, +\, D \, \Delta \left[ u^3 - u  
 -\lambda \, \Delta \, u \right]  & = & 0, \quad \text{ in } \, \Omega_T, 
 \\[0.05in]
 \qquad \ds u &  = & u_0, \quad \text{ on } \,  \partial \Omega_T  , 
 \\
 \textcolor{black}{ \qquad \ds u} &  \textcolor{black}{=} & \textcolor{black}{u_{initial}, \quad \text{ on } \,  \Omega} , 
 \\
 \qquad u^3 - u -\lambda \, \Delta \, u & = & g, \quad \text{ on } \,  \partial \Omega_T, \\
 \end{array}
 \end{array}
}
\end{equation}
where $\Omega_T = \Omega\times(0,T)$ is the space-time domain, \textcolor{black}{ $\partial \Omega_T$ the 
space-time boundary excluding the initial and final time surfaces},
  $D\in\SLOinf$, and 
$\lambda \in\SLOinf$. The diffusion coefficient $D$ and the square width of the transition region $\lambda$ are considered to be constant throughout the domain.
To derive the weak formulation, we use a regular partition 
$\PhT$ of $\Omega_T$ into elements $K_m$, such that:
\begin{equation}
\notag
\label{eq:STdomain}
  \Omega_T = \text{int} ( \bigcup_{\KepT} \overline{K_m} ).
\end{equation}

We apply a mixed FE methodology and introduce two flux variables $\rr,\ttt$ and an
additional scalar variable $q$ as auxiliary variables:
\begin{itemize}
\item $\rr =\{r_x,r_y\}^T= \Nabla u$.
\item $q = u^3 -u -\lambda \, \Nabla \cdot \rr$.
\item$\ttt =\{t_x,t_y\}^T= \Nabla q$.
\end{itemize}
Where $\Nabla$ denotes 
the spatial gradient operator.
Note that the flux variables \textcolor{black}{vary in time due to the definitions of the scalar variables but are only of dimension $\Omega$. }
Hence,  this dictates that the regularity of these trial functions is $\rr\in\SHdivO$, $\ttt\in\SHdivO$, $q\in\SHOOT$, $u\in\SHOOT$, and the IBVP~\eqref{eq:gen_Cahn-Hilliard_BVP} can be recast as an equivalent first-order system of PDEs:
\begin{equation} \label{eq:CH_IBVP_first_order}
\boxed{
\begin{array}{l}
\text{Find }  (u,q,\rr,\ttt) \in \SHOOT\times\SHOOT\times\SHdivO\times\SHdivO \text{ such that:}    
\\[0.05in] 
\qquad 
\begin{array}{rcl}
\ds   \Nabla u - \rr & =  & \bfm{0}, \quad \text{ in } \, \Omega, 
  \\
  \ds   \Nabla q - \ttt & =  & \bfm{0}, \quad \text{ in } \, \Omega, 
  \\
  \ds  u^3 - u - \lambda \, \Nabla \cdot \rr - q & =  & 0, \quad \text{ in } \, \Omega_T, 
  \\
  \ds - \frac{\partial u}{\partial t} +  D \, \Nabla \cdot \ttt & =  & 0, \quad \text{ in } \, \Omega_T, 
 \\[0.025in]
 \qquad u &  = & u_0, \quad \text{ on } \, \partial \Omega_T, 
 \\
 \textcolor{black}{ \qquad \ds u} & \textcolor{black}{ =} & \textcolor{black}{u_{initial}, \quad \text{ on } \,  \Omega} , 
 \\
 \qquad q & = & g, \quad \text{ on } \,  \partial \Omega_T. 
 \end{array}
 \end{array}
}
\end{equation}
\textcolor{black}{The reason we elect to work with this first-order system structure is to apply the DPG
philosophy to construct the optimal test space while using globally continuous FE approximation spaces
such as Lagrange and Raviart-Thomas polynomials for the trial space without the need for auxiliary 
trace unknowns used in the DPG method. } 
We proceed to enforce the PDEs~\eqref{eq:CH_IBVP_first_order} weakly on each
element $\KepT$ and sum the contributions from all $\KepT$, i.e.,
\begin{equation} \label{eq:weak_IBVP_L2_CH}
\begin{array}{c}
\text{Find } \; (u,q,\rr,\ttt) \in \SHOOT\times\SHOOT\times\SHdivO\times\SHdivO:  
\\[0.1in] 
\ds \summa{\textcolor{black}{\KepT}}{} \int_{K_m} \biggl\{\left[ \Nabla u - \rr\right] \cdot \sss_m \, 
\ds+ \left[ \Nabla q - \ttt \right] \cdot \pp_m \,
\ds+ \left[ u^3 - u - \lambda\, \Nabla \cdot \rr - q   \right] \, v_m \\ \,
\ds+ \left[-\frac{\partial u}{\partial t} +  D \, \Nabla \cdot \ttt  \right] \, w_m \biggr\} \dx   = 0, 
 \\[0.2in] \qquad \hspace{2in} \forall (v_m,w_m, \sss_m, \pp_m) \in \textcolor{black}{\SLTOT}\times\textcolor{black}{\SLTOT}\times\SLTOvec\times\SLTOvec.
  \end{array}
\end{equation}
Next, we apply \textcolor{black}{integration by parts} to the terms multiplied with the scalar valued test functions $v_m$ and $w_m$ which dictates 
that we increase the regularity of each scalar valued test function to be in $H^1$ locally for every
$\textcolor{black}{\KepT}$, i.e.,
\begin{equation} \label{eq:weak_IBVP_sum_1_CH}
\begin{array}{c}
\text{Find } \; (u,q,\rr,\ttt) \in \SHOOT\times\SHOOT\times\SHdivO\times\SHdivO:  
\\[0.1in] 
\ds \summa{\textcolor{black}{\KepT}}{}\biggl\{ \int_{K_m}\biggl[ \ds \, \left[ \Nabla u - \rr\right] \cdot \sss_m \, \ds+ \left[ \Nabla q - \ttt \right] \cdot \pp_m \, + \left[ u^3 - u - q  \right] \, v_m \\
+ \ds  \, \lambda \,   \rr \cdot \Nabla v_m
- \ds \frac{\partial u}{\partial t} \, w_m   -  D \, \ttt \cdot \Nabla  w_m \biggr] \dx \\  \ds  \qquad + \oint_{\dKm} D \,  \gamma^m_\nn(\ttt) \, \gamma^m_0(w_m) - \lambda \,  \gamma^m_\nn(\rr) \, \gamma^m_0(v_m) \, \dss \biggr\}  =  0, 
 \\[0.25in]
 \hspace{2in} \forall (v_m,w_m, \sss_m, \pp_m) \in \SHOPT\times\SHOPT\times\SLTOvec\times\SLTOvec,
 \end{array}
\end{equation}
where the broken $H^1$ space on the partition $\PhT$ is defined:
\begin{equation}
\label{eq:broken_h1_space}
\SHOPT \isdef \biggl\{ v\in\SLTOT: \quad v_m \in \SHOK, \; \forall \textcolor{black}{\KepT}\biggr\}.
\end{equation}
The operators $\gamma^m_0: \SHOK: \longrightarrow \SHHdK$ and $\gamma^m_\nn:\SHdivK \longrightarrow \SHmHdK$ denote the trace and normal trace operators (e.g., see~\cite{Girault1986}) on $K_m$; and $\nn_m$ is the outward unit normal vector to the element boundary $\dKm$ of $K_m$. Note that the edge integrals on $\dKm$ are to be interpreted as the duality pairings $\duality{\cdot}{\cdot}{H^{-1/2}(\dKm)}{H^{1/2}(\dKm)}$, instead, we use an integral representation here, as is engineering convention. 

Note that the edge integrals in~\eqref{eq:weak_IBVP_sum_1_CH} only concern the auxiliary flux unknowns $\rr$ and $\ttt$. Thus, any Dirichlet boundary conditions on $u$ and $q$ 
must be enforced strongly. Alternatively, we could perform 
further applications of \textcolor{black}{integration by parts} to shift all the derivatives to the test 
functions, which would allow the weak enforcement of both BCs in~\eqref{eq:CH_IBVP_first_order}.
Finally, these boundary conditions are incorporated in the space $\UUUT$ and we arrive at the AVS-FE weak statement 
for the Cahn-Hilliard IBVP:
\begin{equation} \label{eq:weak_IBVP_sum_2_CH}
\begin{array}{c}
\text{Find } \; (u,q,\rr,\ttt) \in \UUUT:  
\\[0.1in] 
\ds \summa{\textcolor{black}{\KepT}}{}\biggl\{ \int_{K_m}\biggl[ \ds \, \left[ \Nabla u - \rr\right] \cdot \sss_m \, \ds+ \left[ \Nabla q - \ttt \right] \cdot \pp_m \, + \left[ u^3 - u - q  \right] \, v_m \\
+ \ds  \, \lambda \,   \rr \cdot \Nabla v_m
- \ds \frac{\partial u}{\partial t} \, w_m   -  D \, \ttt \cdot \Nabla  w_m \biggr] \dx \\  \ds  \qquad + \oint_{\dKm} D \,  \gamma^m_\nn(\ttt) \, \gamma^m_0(w_m) - \lambda \,  \gamma^m_\nn(\rr) \, \gamma^m_0(v_m) \, \dss \biggr\}  =  0, 
 \\[0.2in]
 \hspace{3in}  \forall (v_m,w_m, \sss_m, \pp_m) \in \VVT,
 \end{array}
\end{equation}
where the trial and test spaces $\UUUT$ and $\VVT$ are defined: 
\begin{equation}
\label{eq:function_spaces_CH}
\begin{array}{c}
\ds \UUUT \isdef \biggl\{ (u,q,\rr,\ttt) \in \SHOOT\times\SHOOT\times\SHdivO\times\SHdivO: \\ \; \gamma_0^m(u)_{|\dKm\cap\partial \Omega_T} =u_0, \gamma_0^m(q)_{|\dKm\cap\partial \Omega_T} =g, \;  \forall\textcolor{black}{\KepT}\biggr\},
\\[0.15in]
\ds \VVT \isdef  \SHOPT\times\SHOPT\times\SLTOvec\times\SLTOvec,
\end{array}
\end{equation}
with  norms $\norm{\cdot}{\UUUT}:  \UUUT \!\! \longrightarrow \!\! [0,\infty)$ and $\norm{\cdot}{\VVT}: \VVT\! \! \longrightarrow\! \! [0,\infty)$ defined as:
\begin{equation}
\label{eq:broken_norms_CH}
\begin{array}{l}
\ds \norm{(u,q,\rr,\ttt)}{\UUUT} \isdef \sqrt{\int_{\Omega_T} \biggl[ (u,u)_{\SHOOT} + (q,q)_{\SHOOT}  +  (\rr,\rr)_{\SHdivO}  +  (\ttt,\ttt)_{\SHdivO}  \biggr] \dx },
\\[0.2in]
\ds   \norm{(v,w, \sss, \pp)}{\VVT} \isdef \\ \qquad \qquad \ds \sqrt{\summa{\textcolor{black}{\KepT}}{}\int_{K_m} \biggl[  h_m^2 \Nabla v_m \cdot \Nabla v_m + v_m^2 + h_m^2 \Nabla w_m \cdot \Nabla w_m + w_m^2    + \sss_m \cdot \sss_m + \pp_m \cdot \pp_m\biggr] \dx }.
 \end{array}
\end{equation}
\textcolor{black}{Note that the scaled norm $\norm{\cdot}{\VVT}$ is equivalent to the $L^2$ 
norm on $\VVT$ (with mesh-dependent equivalence constants):}
\begin{equation}
\label{eq:broken_norm_V}
\begin{array}{l}
\ds  \textcolor{black}{ \norm{(v,w, \sss, \pp)}{V} \isdef } \ds\textcolor{black}{  \sqrt{\summa{\textcolor{black}{\KepT}}{}\int_{K_m} \biggl[   v_m^2  + w_m^2    + \sss_m \cdot \sss_m + \pp_m \cdot \pp_m\biggr] \dx }. }
 \end{array}
\end{equation}
\textcolor{black}{Our choice for the test space norm $\norm{\cdot}{\VVT}$ is motivated by the wish 
to keep all terms in the integral that defines the norm of similar magnitude. This in turn leads to a 
stiffness matrix in which the entries are of similar magnitude which is beneficial for the conditioning 
of the linear system of equations.}
By introducing the operator $B:\UUUT\times\VVT\longrightarrow \mathbb{R}$:
\begin{equation} \label{eq:B_and_F_CH}
\begin{array}{c}
B((u,q,\rr,\ttt);(v,w, \sss, \pp)) \isdef 
\ds \summa{\textcolor{black}{\KepT}}{}\biggl\{ \int_{K_m}\biggl[  \, \left[ \Nabla u - \rr\right] \cdot \sss_m \, + \left[ \Nabla q - \ttt \right] \cdot \pp_m \, + \left[ u^3 - u - q  \right] \, v_m \\
+ \ds  \, \lambda \,   \rr \cdot \Nabla v_m 
- \ds \frac{\partial u}{\partial t} \, w_m -  D \, \ttt \cdot \Nabla  w_m \biggr] \dx
 \\[0.1in]
 \hspace{1in} \biggl. \ds + \oint_{\dKm} D \,  \gamma^m_\nn(\ttt) \, \gamma^m_0(w_m) - \lambda \,  \gamma^m_\nn(\rr) \, \gamma^m_0(v_m) \, \dss \biggr\},
 \end{array}
\end{equation}
we can write the weak formulation~\eqref{eq:weak_IBVP_sum_2_CH} compactly:
\begin{equation} \label{eq:weak_form_CH}
\boxed{
\begin{array}{ll}
\text{Find } (u,q,\rr,\ttt) \in \UUUT   \text{ such that:}
\\[0.05in]
   \quad B((u,q,\rr,\ttt);(v,w, \sss, \pp)) = 0, \quad \forall (v,w, \sss, \pp)\in \VVT. 
 \end{array}}
\end{equation}
\textcolor{black}{See Appendix A for a well-posedness analysis of an AVS-FE weak formulation~\eqref{eq:weak_form_CH} for a linearized Cahn-Hilliard BVP.
}

\subsection{AVS-FE Discretizations}  
\label{sec:avs-discretization}
We seek numerical approximations $(u^h,q^h,\rr^h,\ttt^h)$ of
$(u,q,\rr,\ttt)$ by using FE trial basis functions such as Lagrange 
interpolants or Raviart-Thomas polynomials. 
 However, the test space is 
discretized by employing the DPG philosophy and we use optimal test functions as 
computed from the discrete Riesz problems (see~\eqref{eq:abstract_weak_form_AVS}).
Thus, the FE discretization of~\eqref{eq:weak_form_CH} governing $(u^h,q^h,\rr^h,\ttt^h) \in \UUUhT$ is:
\begin{equation} \label{eq:discrete_form_CH}
\boxed{
\begin{array}{ll}
\text{Find } (u^h,q^h,\rr^h,\ttt^h) \in \UUUhT   \text{ such that:}
\\[0.05in]
   \quad B((u^h,q^h,\rr^h,\ttt^h);(v^h,w^h, \sss^h, \pp^h)) = 0, \quad \forall (v^h,w^h, \sss^h, \pp^h)\in \VVhT, 
 \end{array}}
\end{equation}
where the finite dimensional test space $\VVhT \subset \VVT$ is spanned by numerical
approximations of the test functions through the Riesz representation problems, analogous to~\eqref{eq:local_test_problems}. 

By exploiting the discrete stability of the AVS-FE method, 
the entire space-time domain is discretized by finite elements instead of using traditional time stepping techniques satisfying a CFL condition.
Hence, we have significant flexibility in the choice of mesh parameters in the FE discretization.

To solve the nonlinear variational problem, we can linearize the weak form and solve 
a sequence of linear discrete problems that converge to the nonlinear 
solution. This can be achieved by employing solution procedures such as 
Newton iterations to~\eqref{eq:discrete_form_CH} to which we compute on-the-fly optimal 
test functions at each step of the Newton iterations.
However, we consider an equivalent mixed or saddle point problem interpretation of the AVS-FE method,
as introduced in Section~\ref{sec:AVS-review},
 in which we seek both $(u^h,q^h,\rr^h,\ttt^h)$
and the error representation function $(\psi^h, \varphi^h, \bm{\xi}^h, \bm{\eta}^h)$: 
\begin{equation} \label{eq:disc_mixed}
\boxed{
\begin{array}{ll}
\text{Find } (u^h,q^h,\rr^h,\ttt^h) \in  \UUUhT, (\psi^h, \varphi^h, \bm{\xi}^h, \bm{\eta}^h) \in V^h(\PhT)  \text{ such that:}
\\[0.05in]
   \quad \left(\, (\psi^h, \varphi^h, \bm{\xi}^h, \bm{\eta}^h),(v^h,w^h, \sss^h, \pp^h) \, \right)_\VVT - B((u^h,q^h,\rr^h,\ttt^h);(v^h,w^h, \sss^h, \pp^h)) & =  0, \\& \quad  \hspace*{-2in}\forall (v^h,w^h, \sss^h, \pp^h) \in V^h(\PhT),  \\
  \quad B'_{\mathbbm{u}}((a^h,b^h,\cc^h,\dd^h);(\psi^h, \varphi^h, \bm{\xi}^h, \bm{\eta}^h)) & =  0, \\& \quad \hspace*{-2in}\forall \, ((a^h,b^h,\cc^h,\dd^h)) \in \UUUhT.
 \end{array}}
\end{equation}
Where the operator $B'_{\mathbbm{u}}:\UUUT\times\VVT\longrightarrow \mathbb{R}$ is the 
first order Gateaux derivative of the sesquilinear form $B$ with respect to  
$\mathbbm{u} = (u,q,\rr,\ttt)$.  
Application of the definition of the Gateaux derivative then gives:
\begin{equation} \label{eq:B_prime_CH}
\begin{array}{c}
B'_{\mathbbm{u}}((a,b, \cc,\dd);(\psi, \varphi, \bm{\xi}, \bm{\eta} )) \isdef 
\ds \summa{\textcolor{black}{\KepT}}{}\biggl\{ \int_{K_m}\biggl[  \, \left[ \Nabla a - \cc\right] \cdot \bm{\xi}_m \, + \left[ \Nabla b - \dd \right] \cdot \bm{\eta}_m \, \\ + \left[ 3u^2  a - a - b  \right] \, \psi_m 
+ \ds  \, \beta \,   \cc \cdot \Nabla \psi_m 
- \ds \frac{\partial a}{\partial t} \, \varphi_m -  D \, \dd \cdot \Nabla \varphi_m \biggr] \dx
 \\[0.1in]
 \hspace{1in} \biggl. \ds + \oint_{\dKm} D \,  \gamma^m_\nn(\dd) \, \gamma^m_0(\varphi_m) - \beta \,  \gamma^m_\nn(\cc) \, \gamma^m_0(\psi_m) \, \dss \biggr\}
 \\[0.15in]
 \end{array}
\end{equation}

\
\section{Numerical Verifications}  
\label{sec:verifications}
In this section, we present several numerical verifications applying the AVS-FE method to 
\textcolor{black}{stationary and transient} problems. 
To establish the  solution of~\eqref{eq:disc_mixed} we use the high-level
FE solvers Firedrake~\cite{rathgeber2017firedrake} and FEniCS~\cite{alnaes2015fenics} which in turn employ the Portable, Extensible Toolkit for Scientific Computation 
(PETSc) library Scalable Nonlinear Equations Solvers (SNES)~\cite{abhyankar2018petsc,petsc-user-ref} to perform Newton iterations.
In all experiments presented, we use PETSc SNES objects in Firedrake~\cite{rathgeber2017firedrake} and
FEniCS~\cite{alnaes2015fenics} with the default settings for tolerances for the iterations. 

\begin{figure}[h]
{\centering
\includegraphics[width=0.5\textwidth]{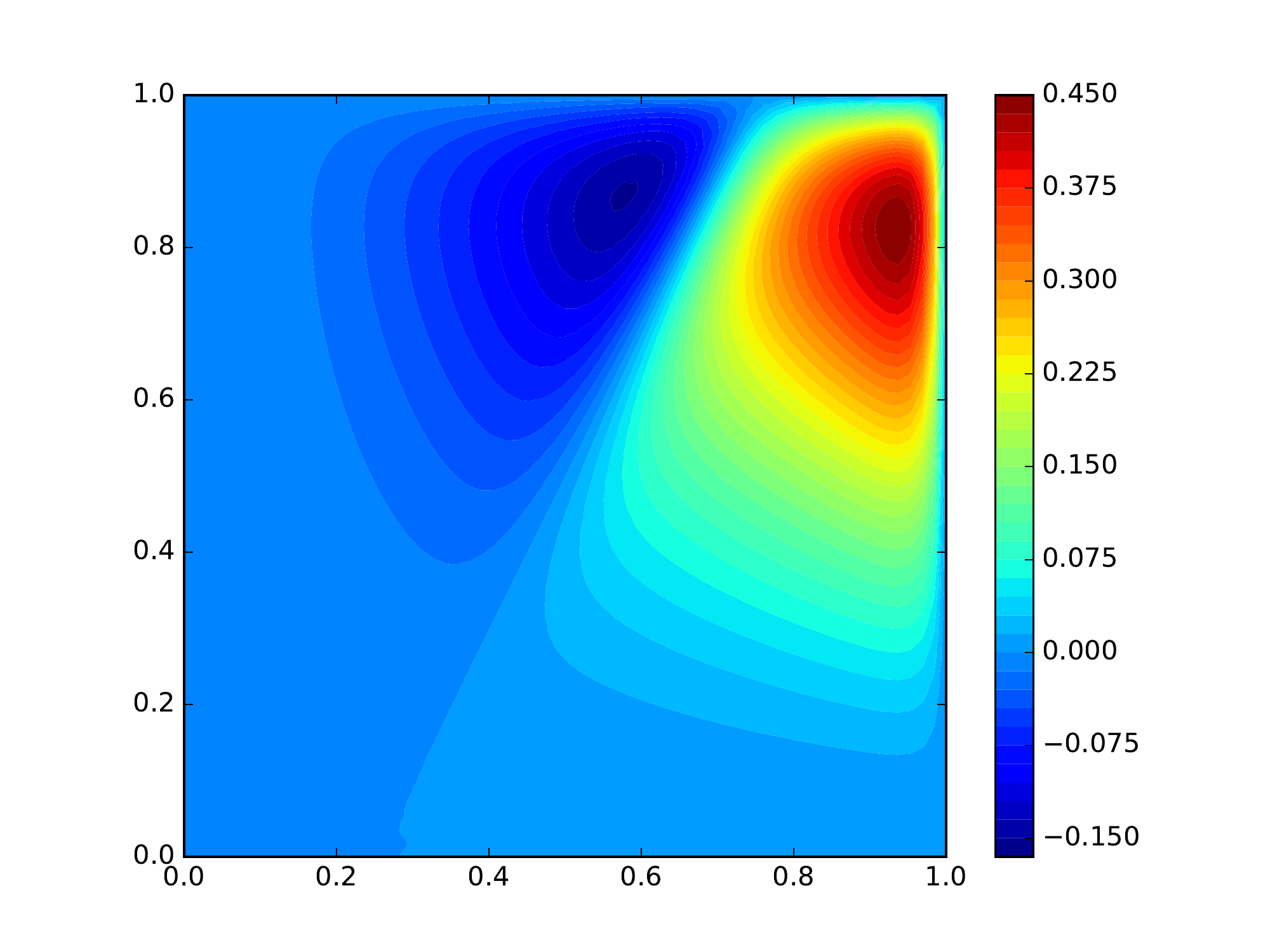}}
 \caption{\label{fig:prop_exact} \textcolor{black}{ Exact solution $u(x,y)$ of the propagating front problem. }}
\end{figure}

\textcolor{black}{
\subsection{Numerical Convergence Studies}
\label{sec:1d_prob}
To ascertain the convergence behavior of the AVS-FE method for the Cahn-Hilliard equation, we perform 
multiple verifications of its convergence properties. 
We first consider a stationary model problem where we consider a manufactured exact solution
$u(x,y)$ slightly modified from~\cite{van2011goal} called the \emph{propagating front test case}}:
\begin{equation}
\label{eq:prop_fron}
\begin{array}{c}
\ds \textcolor{black}{ u(x,y)} = \textcolor{black}{(x \,  y)} \, \text{\textcolor{black}{tanh}} \color{black}{ \left(  \textcolor{black}{\frac{\ds x - 0.5y -0.25}{\ds\sqrt{2 \lambda} } } \right)  \left(x + \frac{e^{50 x}-1}{1-e^{50}}  \right) \left( y + \frac{e^{10 x}-1}{1-e^{10}}  \right).}
\end{array}
\end{equation}
\textcolor{black}{To impose this exact solution, we apply the Cahn-Hilliard equation~\eqref{eq:Cahn-Hilliard_BVP} 
to~\eqref{eq:prop_fron} to ascertain a corresponding nonzero right hand side and Dirichlet boundary
conditions on $u$ and $q$, we pick $D=1$, $\lambda = 1/320$, and the domain as the unit square i.e., $\Omega = (0,1)\times(0,1)$. 
This exact solution is shown in Figure~\ref{fig:prop_exact} }.
\textcolor{black}{ The function spaces we use for this verification consist of Raviart-Thomas and
 Lagrangian bases for the $\SHdivO$ and $\SHOO$ approximations, respectively as well as their 
 discontinuous counterparts for $\SHdivP$ and $\SHOP$. In the mixed problem 
we solve~\eqref{eq:disc_mixed}, we pick basis functions for the solution $(u^h,q^h,\rr^h,\ttt^h)$ and
the error representation function $(\psi^h, \varphi^h, \bm{\xi}^h, \bm{\eta}^h)$ of identical
approximation order. }


\textcolor{black}{ We implement the two-dimensional problem in FEniCS 
and start with a uniform mesh consisting of two triangular elements to 
which we perform both uniform and adaptive $h-$refinements.
For uniform mesh refinements, we establish the corresponding rates of convergence to ensure 
optimal behavior. For the base variable $u$, the expected rates of convergence 
in Sobolev type norms for a linear fourth order PDE are:} 
\begin{equation}
\label{eq:conv_rates}
\boxed{  
\begin{array}{rlr}
\ds \norm{u-u^h}{\SLTO} \leq & C \, h^{p} & \text{ if } p < 3,
 \\[0.1in]
 \ds \norm{u-u^h}{\SLTO} \leq & C \, h^{p+1} & \text{ if } p \geq 3,
 \\[0.1in]
 \ds \norm{u-u^h}{\SHOO} \leq & C \, h^{p}.
 \\[0.1in]
 \end{array} }
\end{equation}
\textcolor{black}{See~\cite{kastner2016isogeometric} and references therein} for error estimates of the Cahn-Hilliard equation \textcolor{black}{and fourth order PDEs}.
Note that for $p < 3$, the $L^2$ and $H^1$ norms of the approximation error converge at the 
same rate \textcolor{black}{since the estimates depend upon the order of the PDE $2m=4$}.

\textcolor{black}{
Due to the implementation of the AVS-FE method as a mixed problem~\eqref{eq:disc_mixed}, we establish both the
 approximate solutions $u^h,q^h,\rr^h,\ttt^h$ and the 
error representation function $(\psi^h, \varphi^h, \bm{\xi}^h, \bm{\eta}^h)$. 
Hence, for the $h-$adaptive algorithm, we use the restriction of this estimate to each element as
 an error indicator:} 
\begin{equation}
\label{eq:err_indicator}
\begin{array}{c}
\ds \textcolor{black}{\eta = \norm{(\psi^h, \varphi^h, \bm{\xi}^h, \bm{\eta}^h)}{\VVK} }
\end{array}
\end{equation}
\textcolor{black}{This error indicator 
has been successfully applied for the AVS-FE method for the linear convection-diffusion PDE
and to several classes of problems of the DPG~\cite{Demkowicz2,Demkowicz6}. 
In particular, in~\cite{Demkowicz2}, (see Theorem 2.1) Carstensen \emph{et al.} investigate and verify the 
robustness of this error estimate and corresponding indicators under the requirement of existence 
of a Fortin operator~\cite{nagaraj2017construction,demkowicz2020construction}.
Without such a Fortin operator, discrete stability and convergence of DPG and AVS-FE methods 
would not be possible. To mark elements for refinement, we consider the strategy of 
D\" orfler~\cite{dorfler1996convergent}, based on the approximate total energy error. 
The stabilized adaptive method of Calo \emph{et al.} introduced in~\cite{calo2019adaptive} also
utilizes this type refinement strategy and error indicator. }

\begin{figure}[h]
\subfigure[ \label{fig:convegence_results_CHp1_L2} $\norm{u-u^h}{\SLTO}$. ]{\centering
 \includegraphics[width=0.44\textwidth]{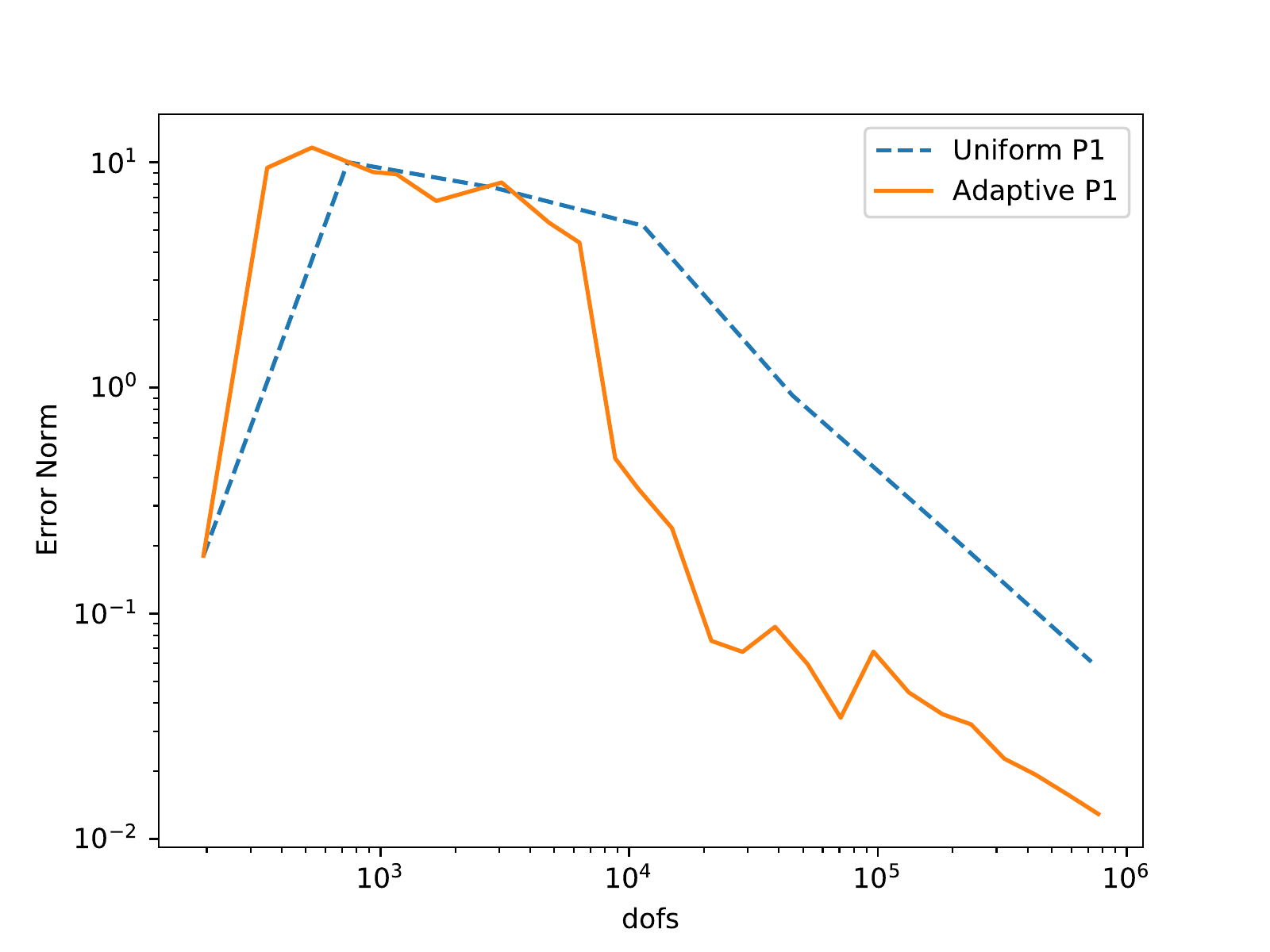}}
  \subfigure[ \label{fig:convegence_results_CHp1_H1}  $\norm{q-q^h}{\SLTO}$.]{\centering 
 \includegraphics[width=0.44\textwidth]{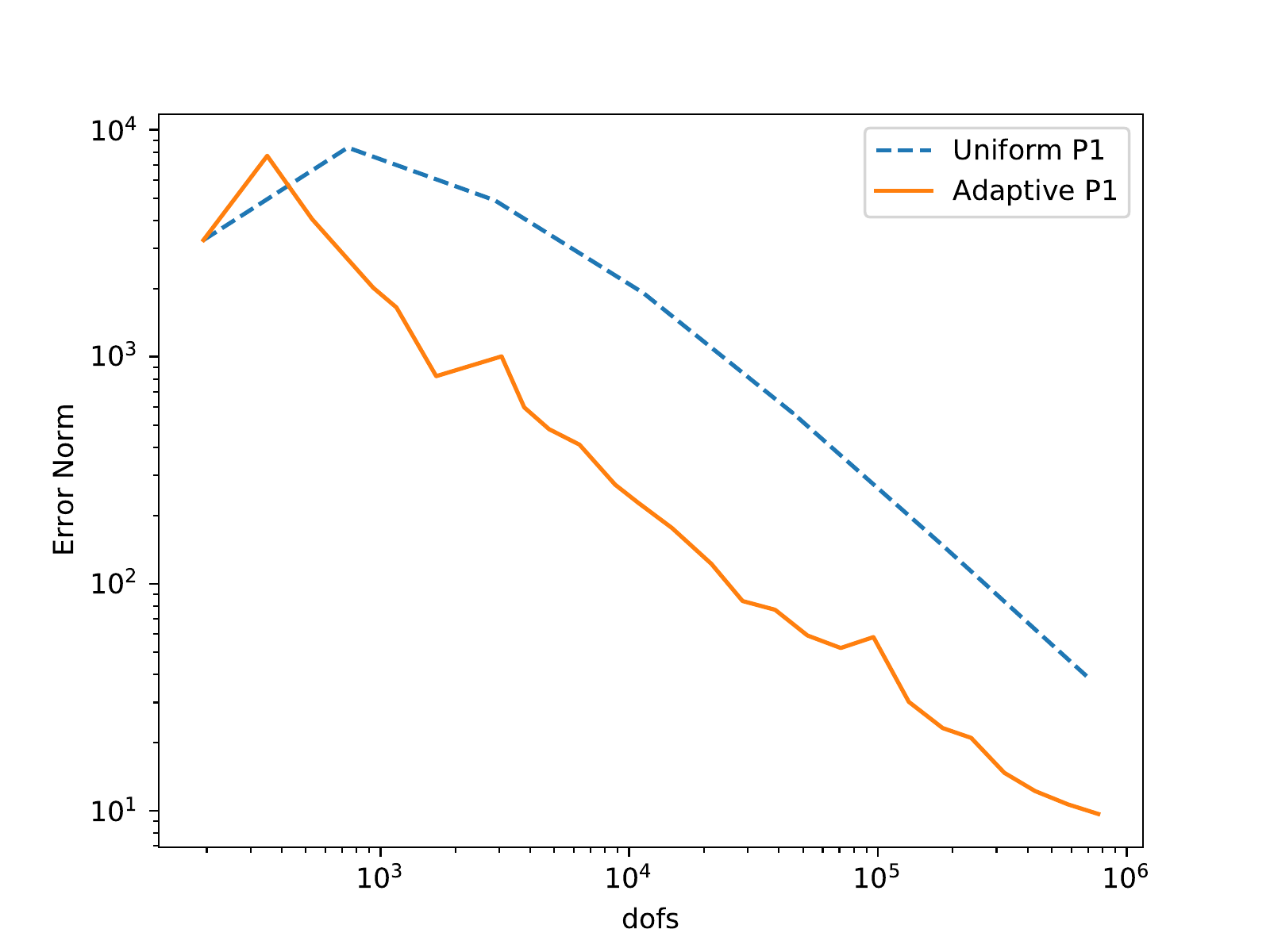}}
  \caption{\label{fig:convegence_results_CHp1} Convergence of uniform and adaptive refinements for the stationary case, linear approximations.}  
\end{figure}
\begin{figure}[h]
\subfigure[ \label{fig:convegence_results_CHp2_L2} $\norm{u-u^h}{\SLTO}$. ]{\centering
 \includegraphics[width=0.44\textwidth]{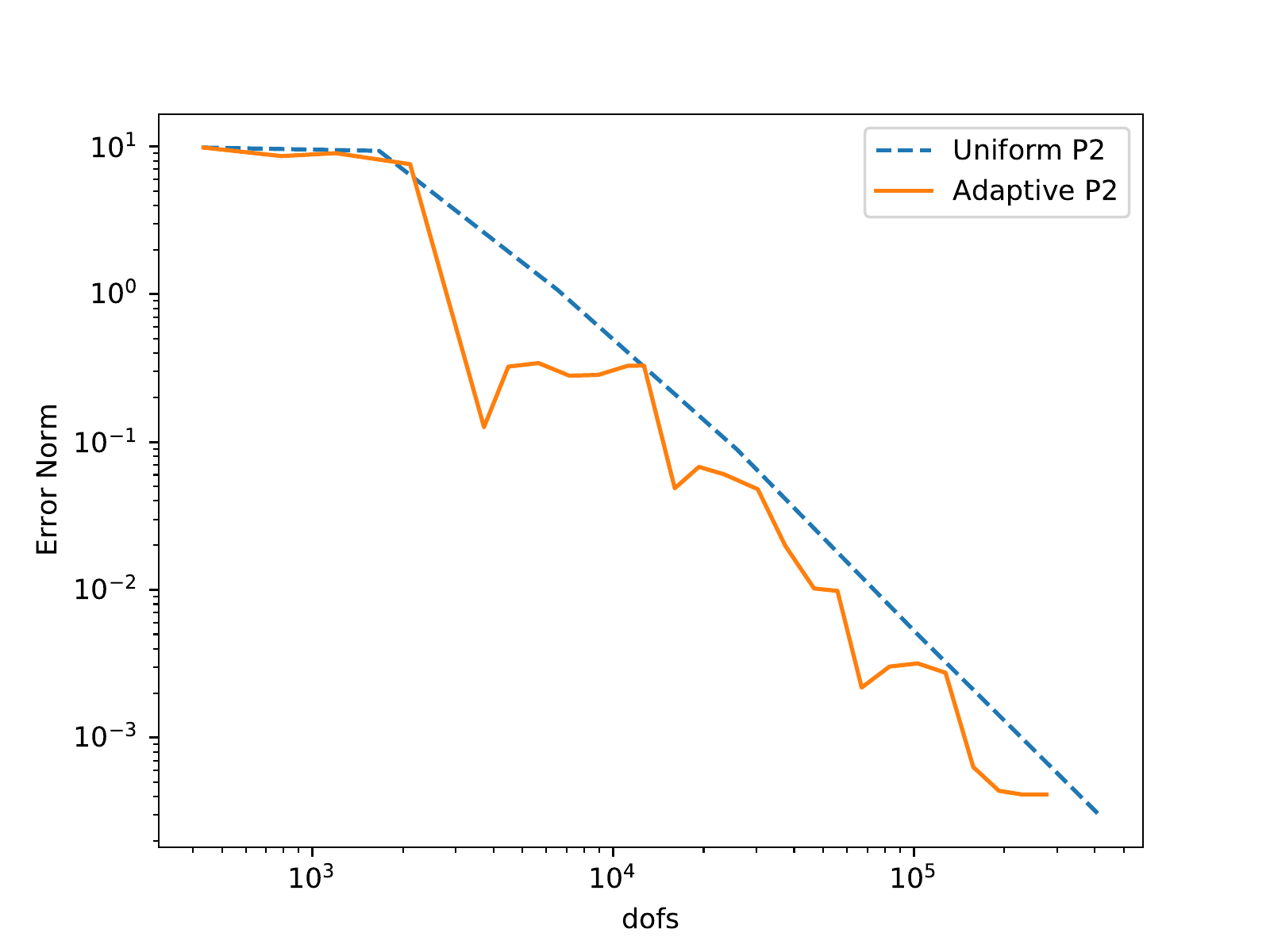}}
  \subfigure[ \label{fig:convegence_results_CHp2_H1}  $\norm{q-q^h}{\SLTO}$.]{\centering 
 \includegraphics[width=0.44\textwidth]{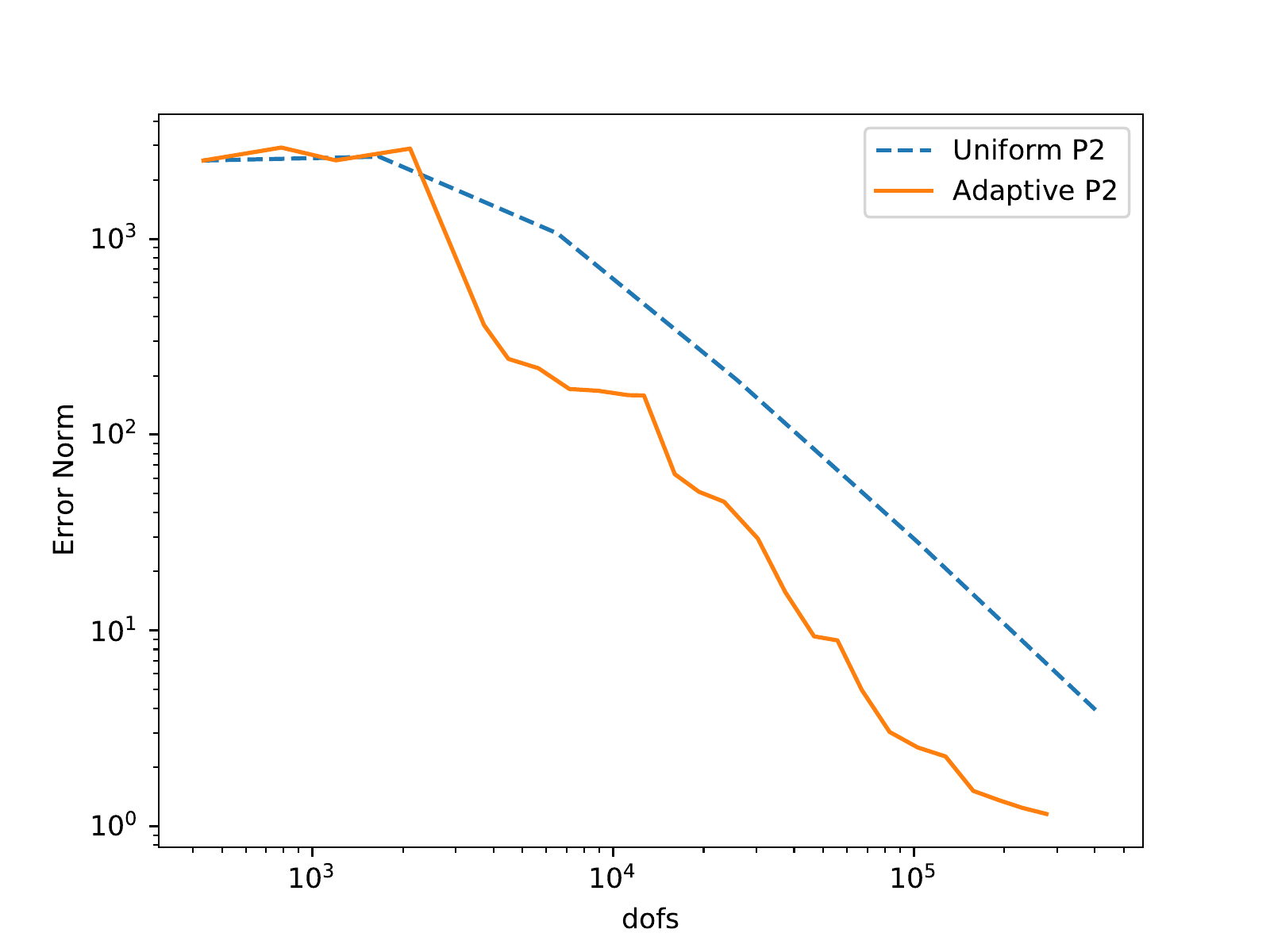}}
  \caption{\label{fig:convegence_results_CHp2} Convergence of uniform and adaptive refinements for the stationary case, quadratic approximations.}  
\end{figure}
\begin{figure}[h]
\subfigure[ \label{fig:convegence_results_CHp3_L2} $\norm{u-u^h}{\SLTO}$. ]{\centering
 \includegraphics[width=0.44\textwidth]{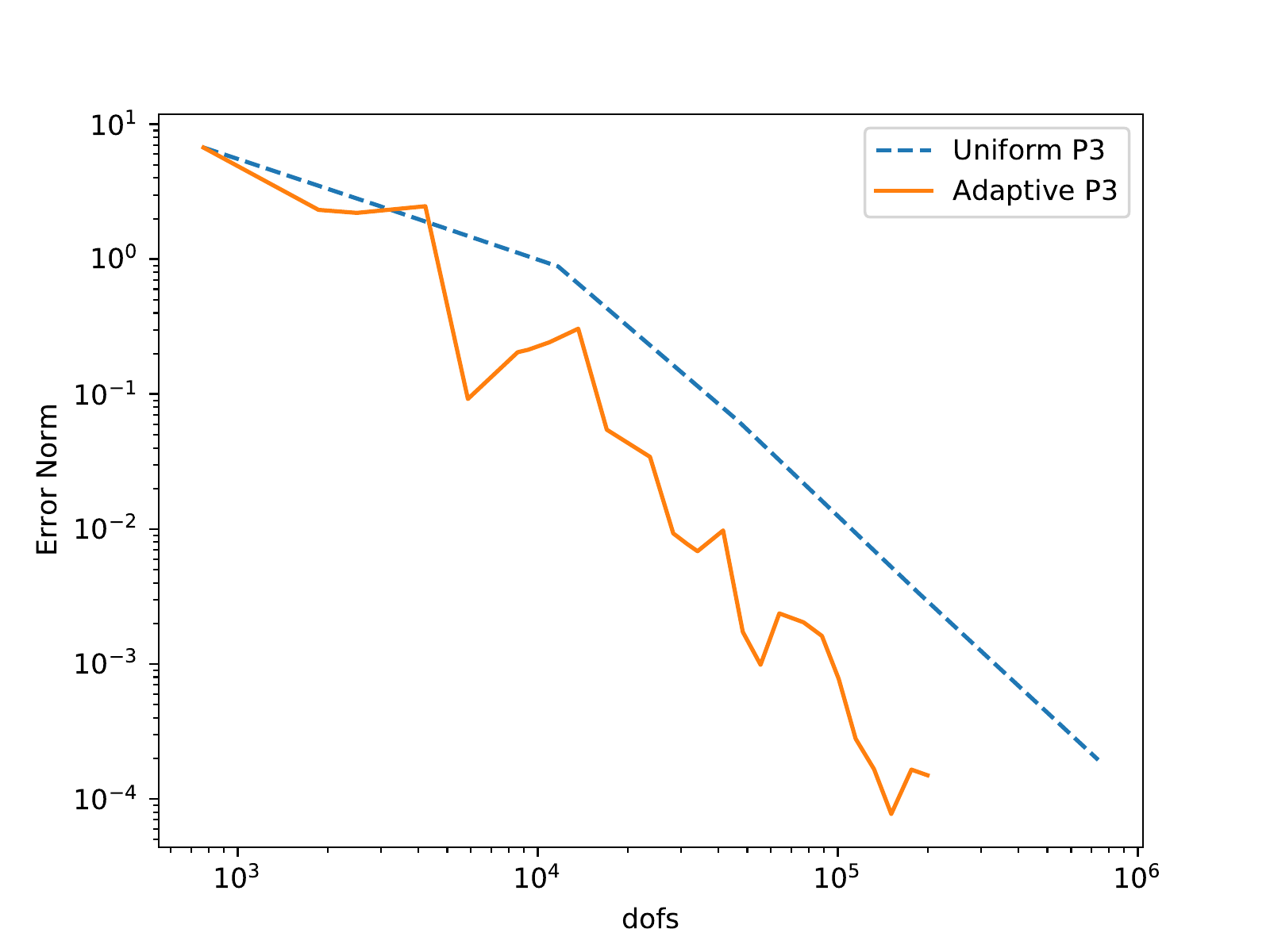}}
  \subfigure[ \label{fig:convegence_results_CHp3_H1}  $\norm{q-q^h}{\SLTO}$.]{\centering 
 \includegraphics[width=0.44\textwidth]{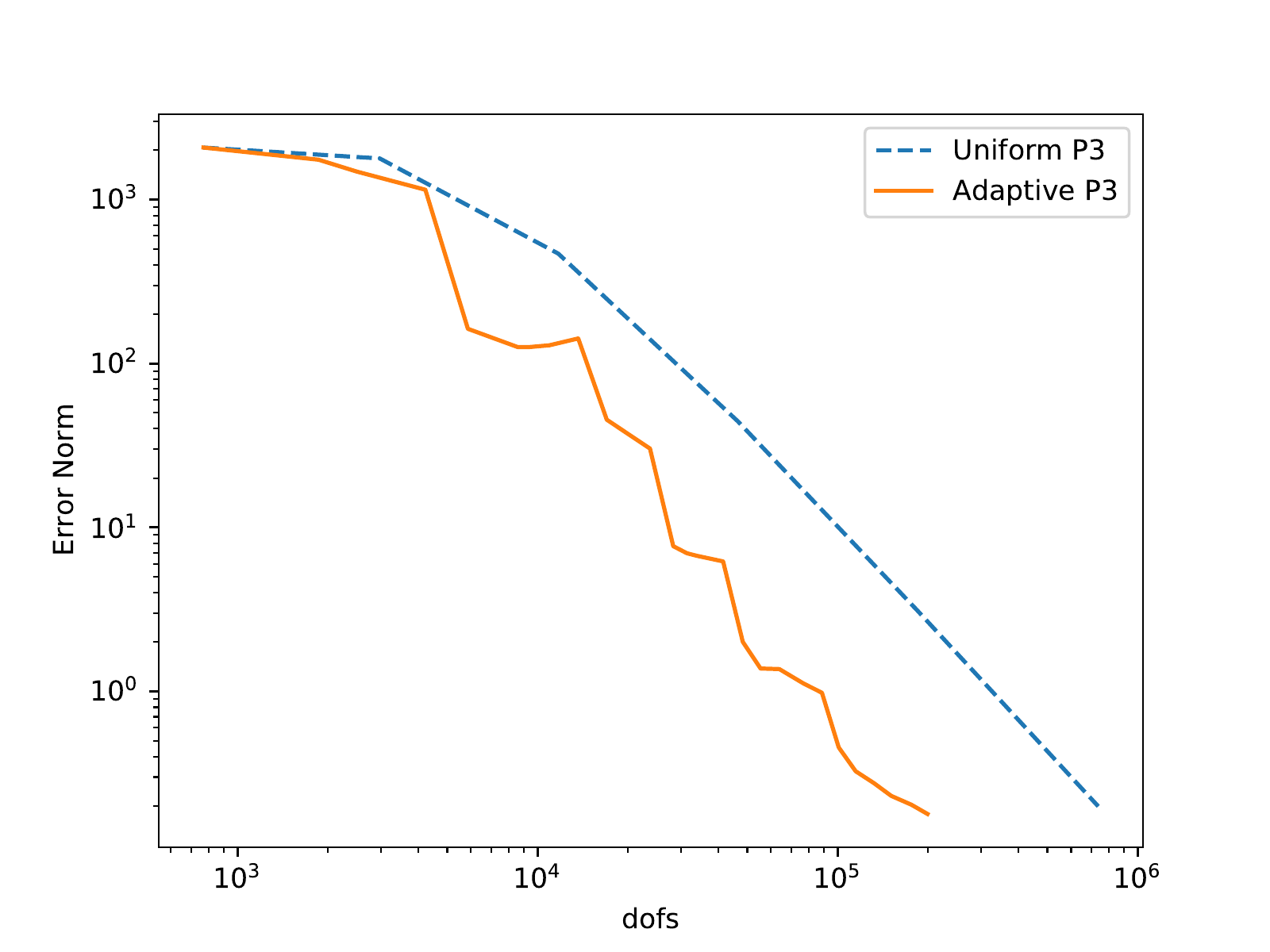}}
  \caption{\label{fig:convegence_results_CHp3} Convergence of uniform and adaptive refinements for the stationary case, cubic approximations.}  
\end{figure}

\textcolor{black}{In Figures~\ref{fig:convegence_results_CHp1},~\ref{fig:convegence_results_CHp2}, and~\ref{fig:convegence_results_CHp3} we present the convergence histories for both adaptive and uniform 
mesh refinements for increasing orders of approximation.}
\textcolor{black}{The observed uniform convergence rates in $\norm{u-u^h}{\SLTO}$ for the uniform refinements are $h^{p+0.9}, h^{p+2}$ and $h^{p+1}$, for linear, quadratic and cubic approximations, respectively.
 For the linear and quadratic 
case these rates are higher than expected in~\eqref{eq:conv_rates}, whereas the convergence 
rates in the $H^1$ norm are identical to those in~\eqref{eq:conv_rates}. 
The observed uniform convergence rates in the $L^2$ error $\norm{q-q^h}{\SLTO}$ for the uniform refinements are $h^{p}, h^{p+1}$ and $h^{p+0.7}$, for linear, quadratic and cubic approximations, respectively.
In all cases, the adaptive refinements leads to lower errors than the uniform refinements as expected as shown in Figures~\ref{fig:convegence_results_CHp1},~\ref{fig:convegence_results_CHp2}, and~\ref{fig:convegence_results_CHp3}. }

\begin{figure}[h]
\subfigure[ \label{fig:comp_energy} $\norm{(u,q,\rr,\ttt)-(u^h,q^h,\rr^h,\ttt^h)}{B}$. ]{\centering
 \includegraphics[width=0.44\textwidth]{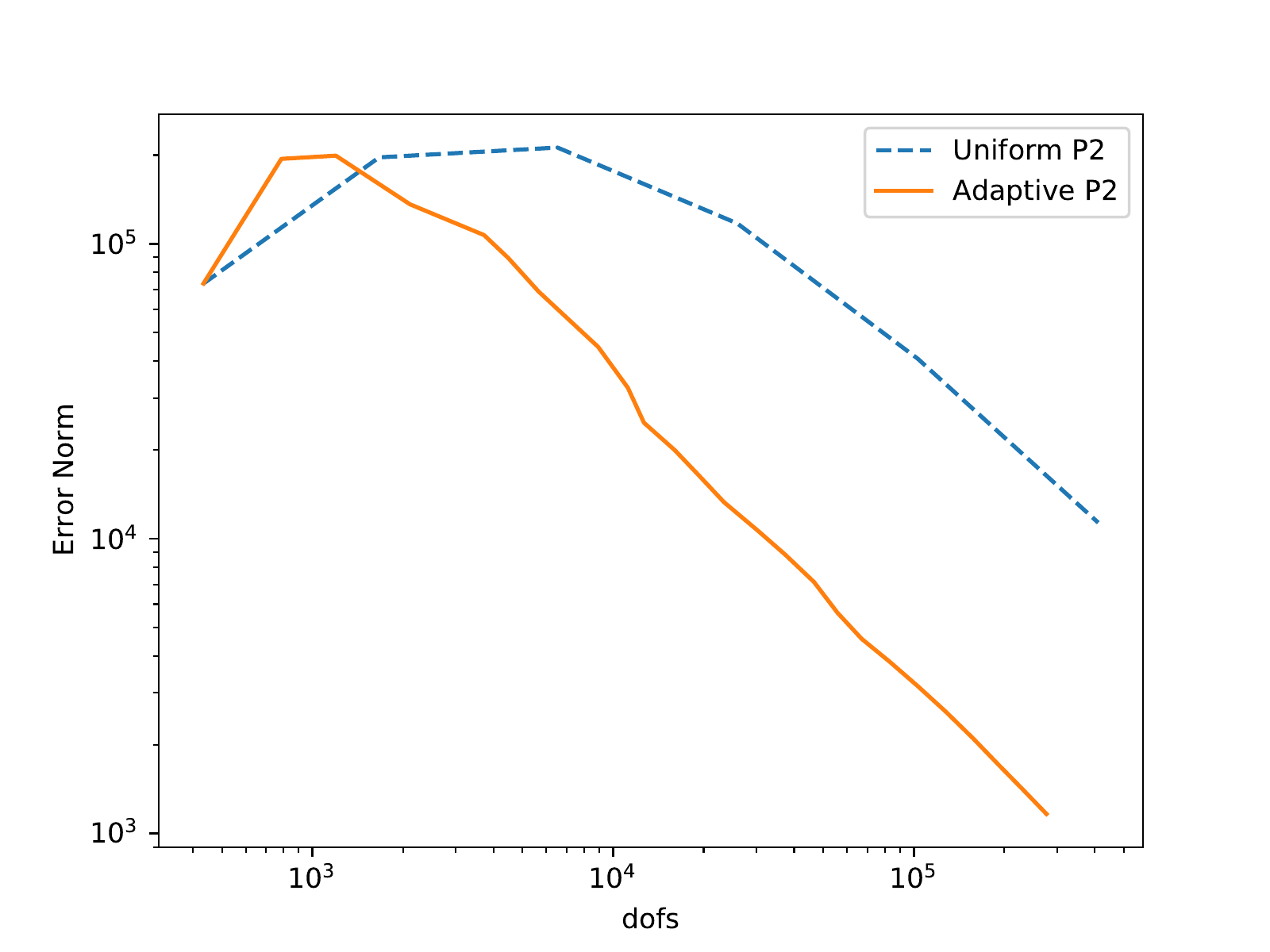}}
  \subfigure[ \label{fig:comp_r_gradu}  Gradient and flux approximation comparison.]{\centering 
 \includegraphics[width=0.44\textwidth]{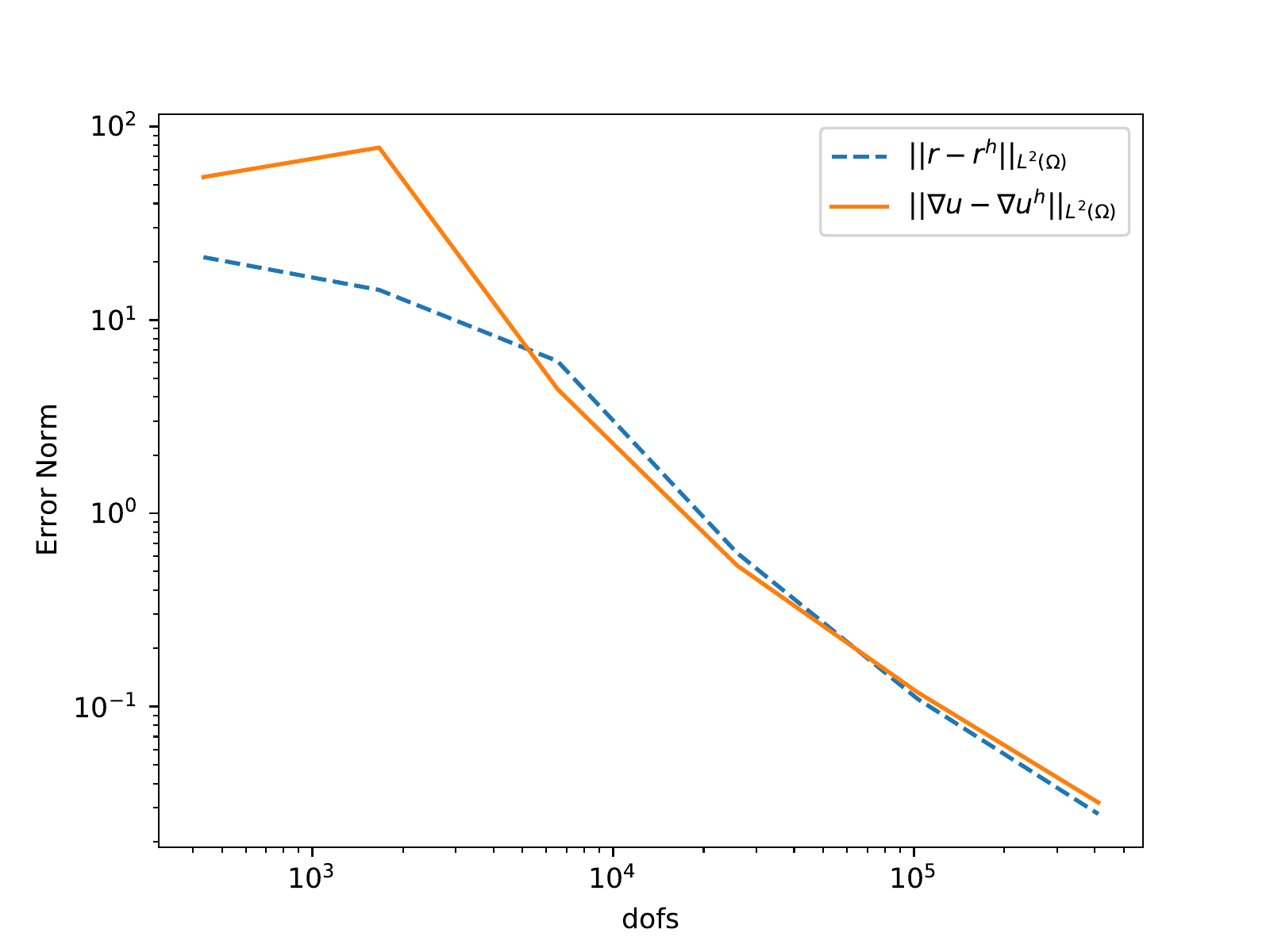}}
  \caption{\label{fig:convegence_results_CH1} Convergence results for the stationary Cahn-Hilliard 
  problem.}
\end{figure}

\textcolor{black}{While the adaptive refinement strategy delivers lower errors than uniform refinements in terms of the 
$L^2$ errors of the base variable $u$, the difference between uniform and adaptive refinements is 
significantly more noticeable in terms of the energy norm as shown in Figure~\ref{fig:comp_energy}.
The difference between the two curves is about an order of magnitude. The reason for this large 
disparity between the two refinement procedures is our choice of error indicator, which is a local 
representation of the energy norm, as well as the refinement criterion based on the total energy 
error. Hence, in the adaptive strategy, the goal is to minimize the energy error.   
In Figure~\ref{fig:comp_r_gradu} we compare the $L^2$ errors of $\Nabla u$ and the vector 
variable $\rr$ for the case of quadratic approximations under uniform refinement. The difference 
between the two is marginal in this case with the vector variable being slightly more accurate for the 
finest mesh.}
%
%
%
%
%

\begin{figure}[h!]
\subfigure[ \label{fig:in_mesh} Mesh. ]{\centering
 \includegraphics[width=0.45\textwidth]{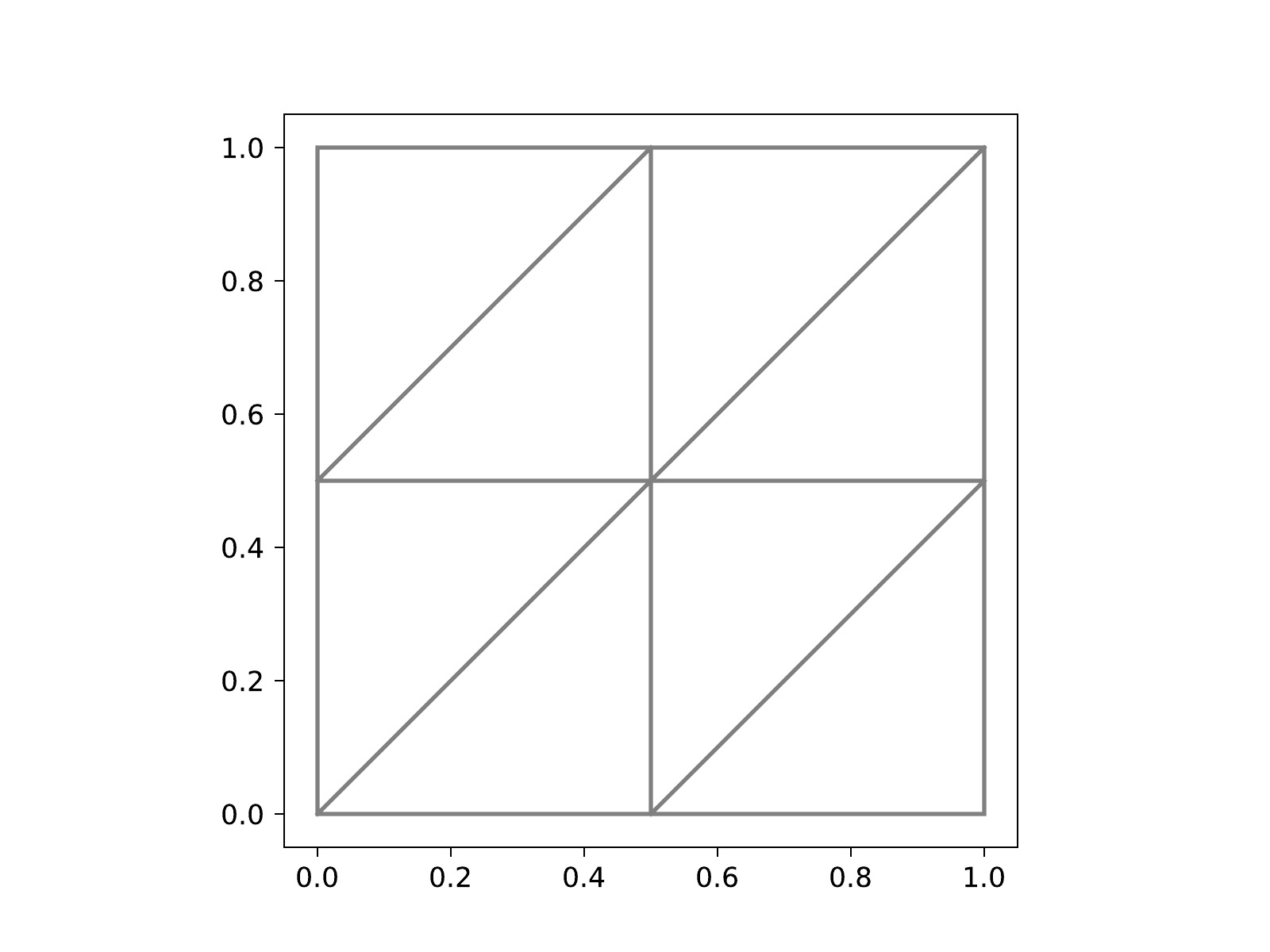}}
  \subfigure[ \label{fig:in_sol} Solution $u^h$.]{\centering
 \includegraphics[width=0.45\textwidth]{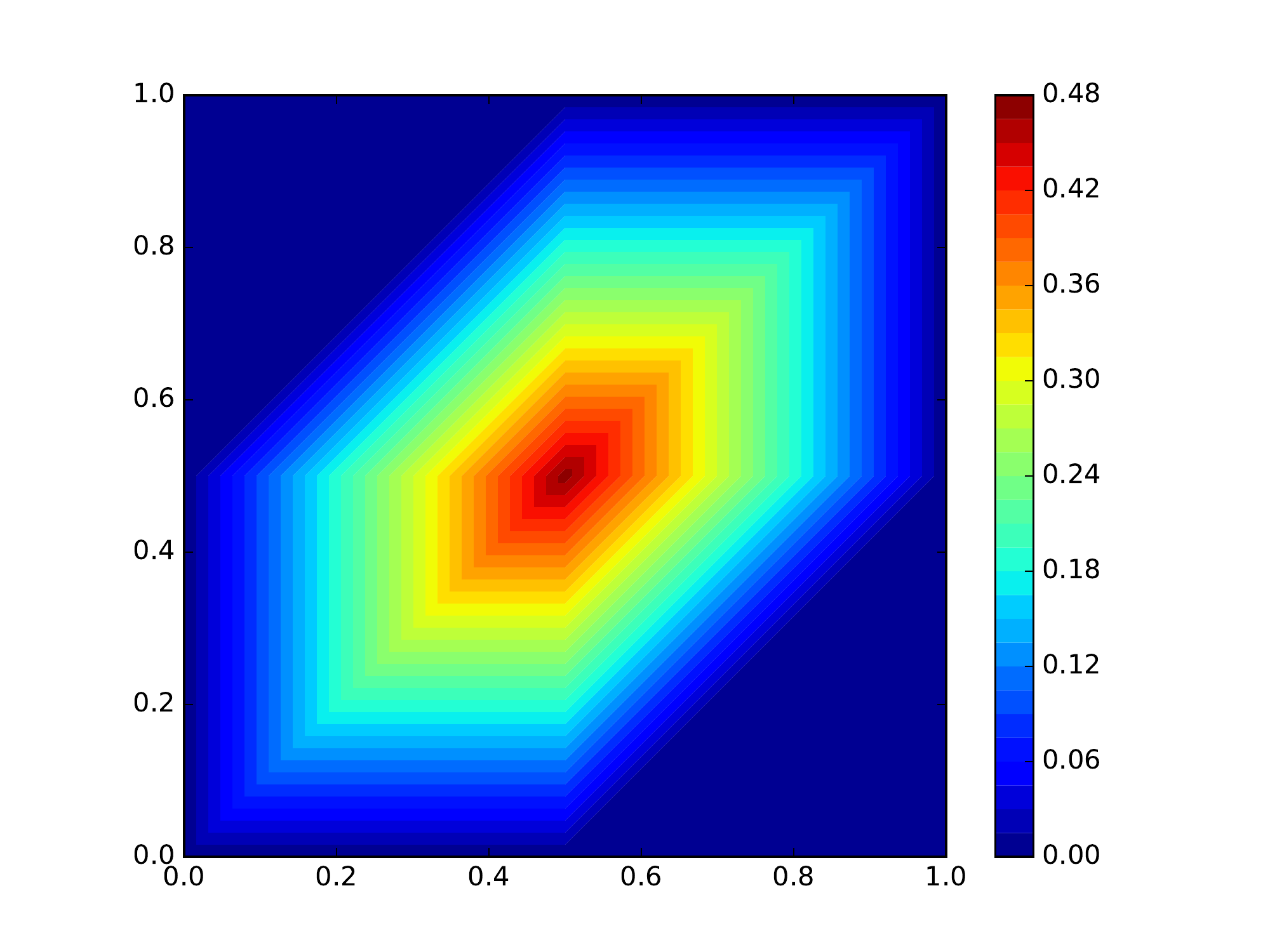}}
  \caption{\label{fig:initial_adaptive_step} Initial adaptive step.}
\end{figure}

\begin{figure}[h!]
\subfigure[ \label{fig:14_mesh} Mesh. ]{\centering
 \includegraphics[width=0.45\textwidth]{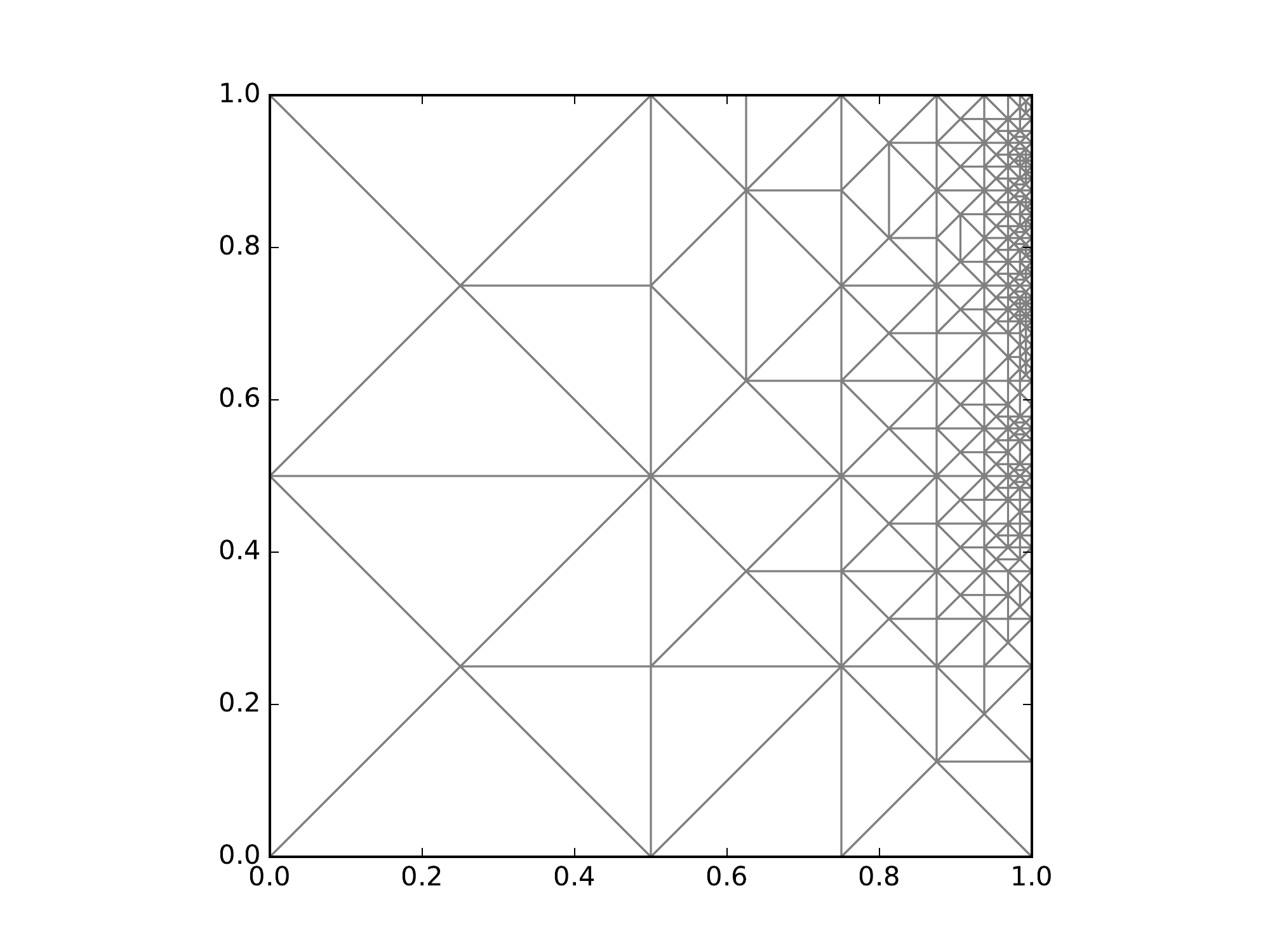}}
  \subfigure[ \label{fig:14_sol} Solution $u^h$.]{\centering
 \includegraphics[width=0.45\textwidth]{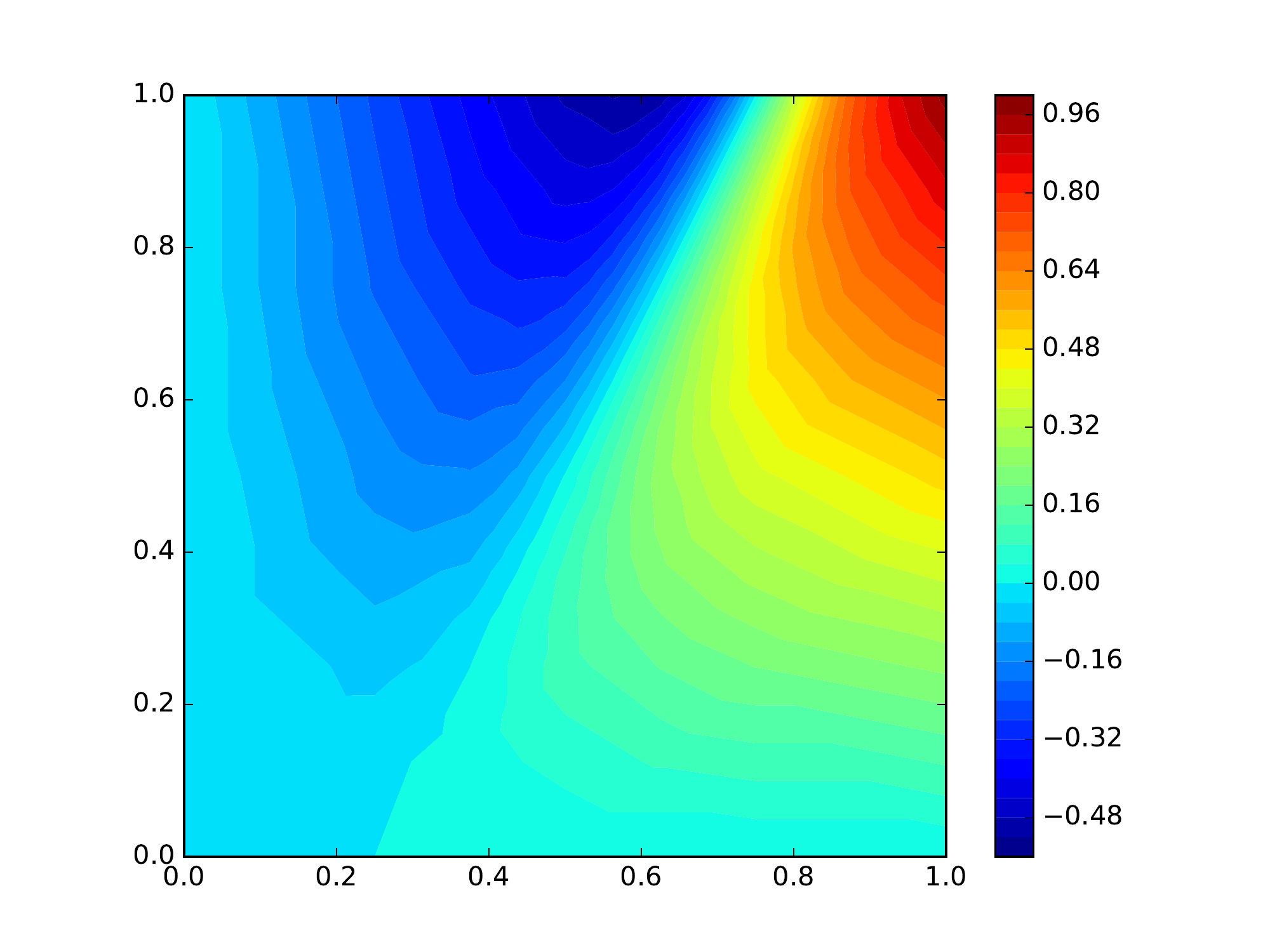}}
  \caption{\label{fig:14_adaptive_step} 10'th adaptive step.}
\end{figure}

\begin{figure}[h!]
\subfigure[ \label{fig:24_mesh} Mesh. ]{\centering
 \includegraphics[width=0.45\textwidth]{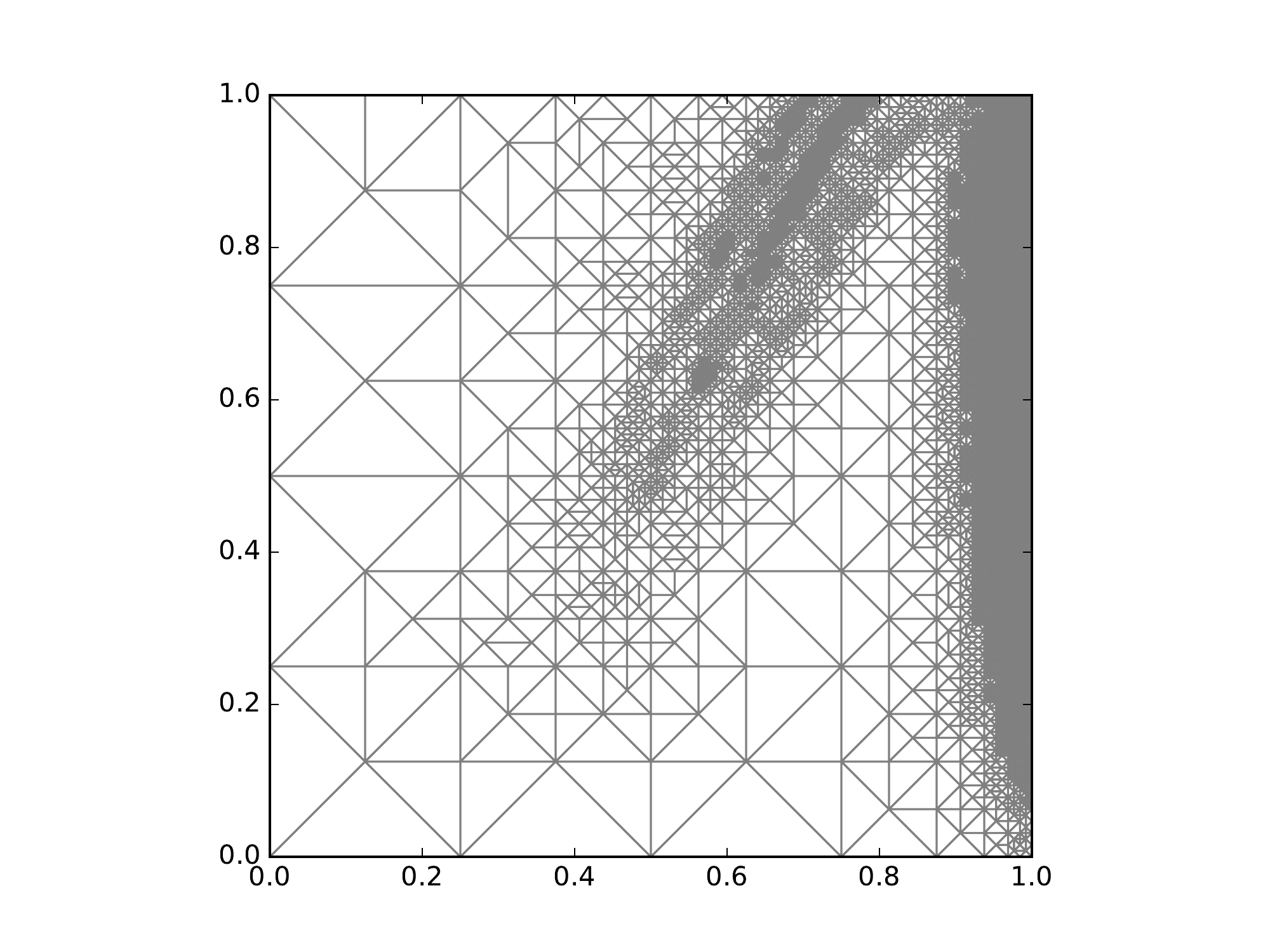}}
  \subfigure[ \label{fig:24_sol} Solution $u^h$.]{\centering
 \includegraphics[width=0.45\textwidth]{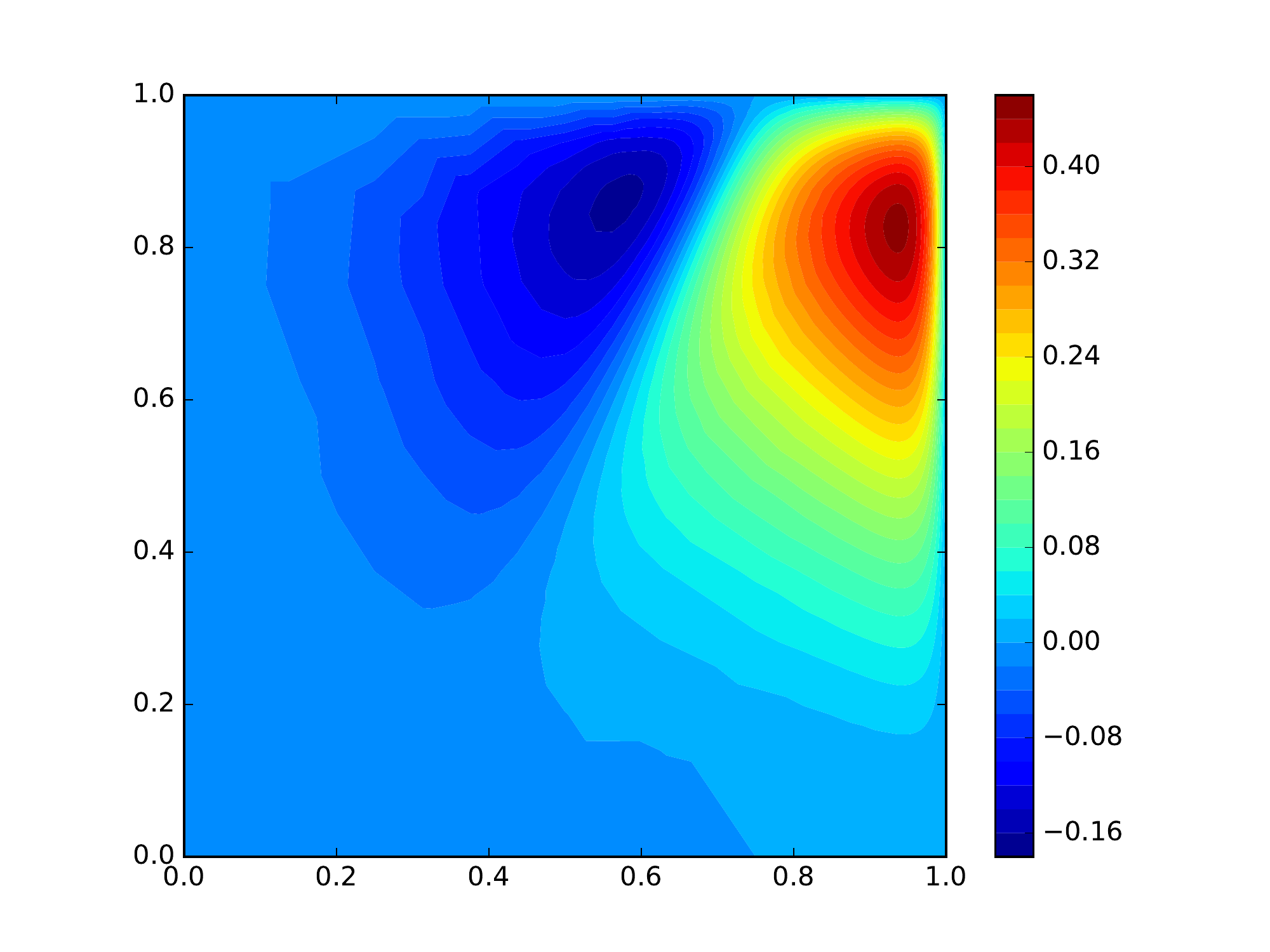}}
  \caption{\label{fig:24_adaptive_step} 25'th adaptive step.}
\end{figure}

\textcolor{black}{To present a sequence of adaptively refined meshes, we consider the case of first order approximations.
In Figures~\ref{fig:initial_adaptive_step},~\ref{fig:14_adaptive_step},
and~\ref{fig:24_adaptive_step} we show selected meshes and 
corresponding solutions from the 
refinement process.} The built-in error indicator performs very well as the mesh refinements are focused along the
propagating front.
%
%




\textcolor{black}{Next, we consider a transient model problem to verify the convergence behavior of 
of the space-time AVS-FE method. To this end, we consider the exact solution:} 
\begin{equation}
\label{eq:2dconv}
\begin{array}{c}
\ds \textcolor{black}{ u(x,y,t) =}     \color{black}{  
\text{sin}(\pi t) \, \text{sin}(\pi y) \, \text{sin}(\pi x) , }
\end{array}
\end{equation}
\textcolor{black}{and consider the space-time domain $\Omega_{T} = (0,1)\times(0,1)\times(0,0.1)$. 
This exact solution is used to ascertain a nonzero source term by applying the Cahn-Hilliard 
differential operator as well as initial and Dirichlet boundary conditions. The parameters $D$ and 
$\lambda$ are both chosen to be unity for simplicity. 
In the same fashion as the stationary example, we compare both uniform and adaptive mesh refinements 
based on the same principles. We discretize all continuous variables using equal order Lagrange  polynomial
functions in this case to show its effect on the accuracy of the flux variables as compared 
to the Raviart-Thomas approximations in the previous case. Since the source resulting from the 
transient exact solution~\eqref{eq:2dconv} is smooth, the resulting regularity of the flux 
variables $\rr,\ttt$ is higher than $H(div)$. Hence, the increased regularity of the approximation 
will not lead to consistency issues.
All components of the error representation function are discretized using discontinuous Lagrange 
polynomials of the same order as the continuous trial variables.}
\begin{figure}[h]
\subfigure[ \label{fig:trans_comp_energy_p1} $\norm{(u,q,\rr,\ttt)-(u^h,q^h,\rr^h,\ttt^h)}{\UUUT}$. ]{\centering
 \includegraphics[width=0.44\textwidth]{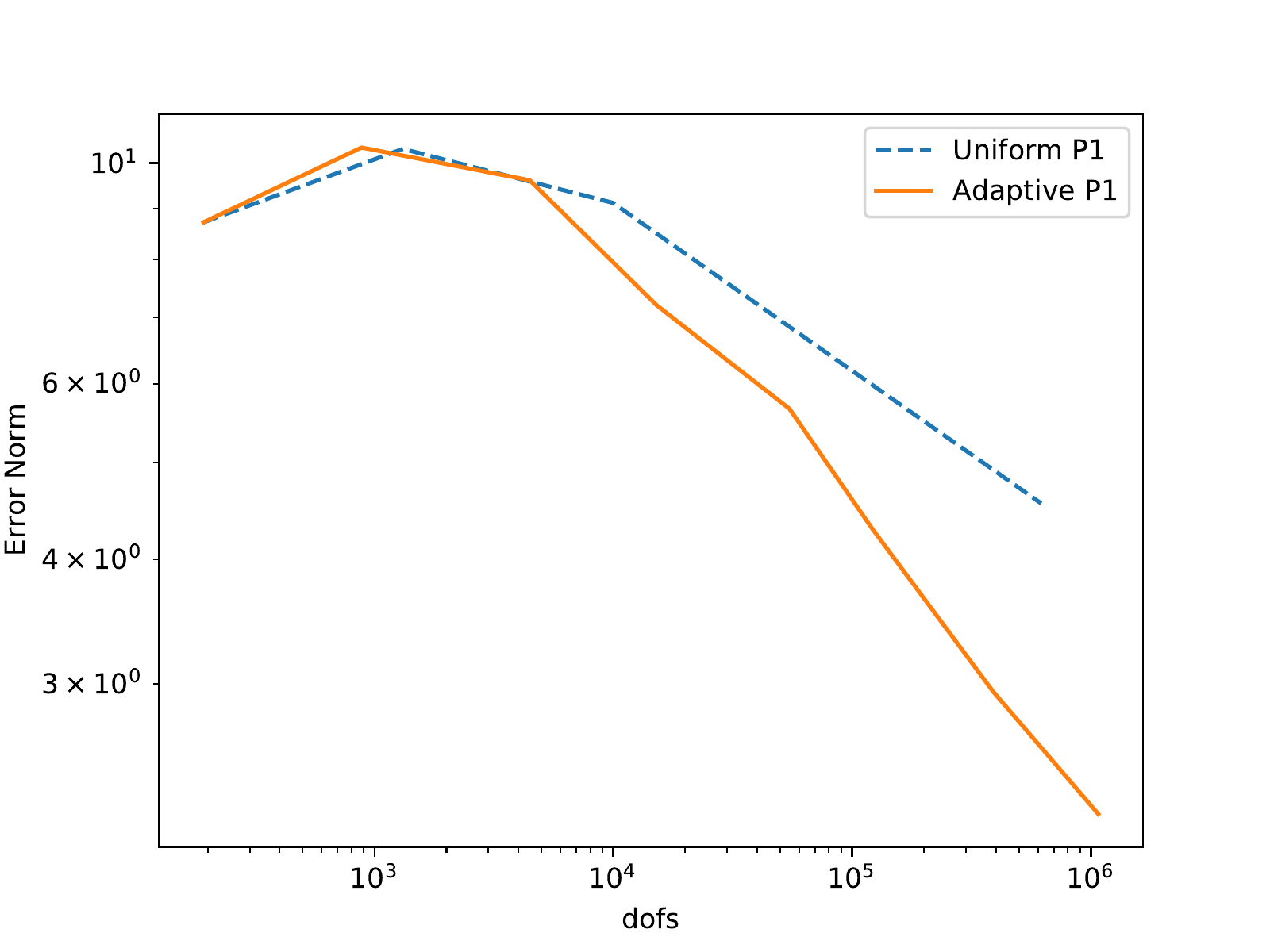}}
  \subfigure[ \label{fig:trans_comp_r_gradu_p1}  Gradient and flux approximation comparison for adaptive refinements.]{\centering 
 \includegraphics[width=0.44\textwidth]{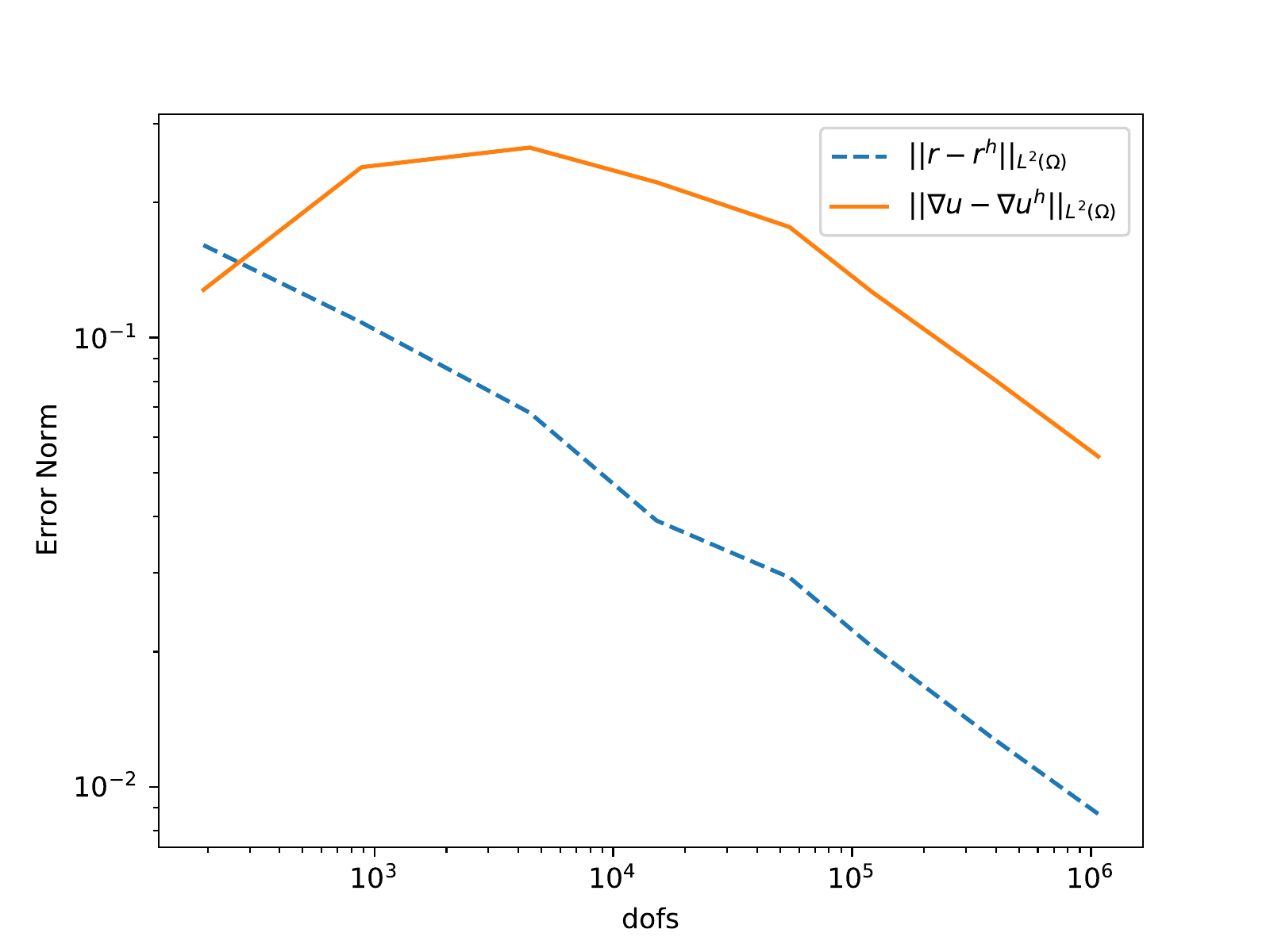}}
  \caption{\label{fig:trans_convegence_results_CH1} Convergence results for the transient Cahn-Hilliard 
  problem for linear basis functions.}
\end{figure}
\begin{figure}[h]
\subfigure[ \label{fig:trans_comp_energy_p2} $\norm{(u,q,\rr,\ttt)-(u^h,q^h,\rr^h,\ttt^h)}{\UUUT}$. ]{\centering
 \includegraphics[width=0.44\textwidth]{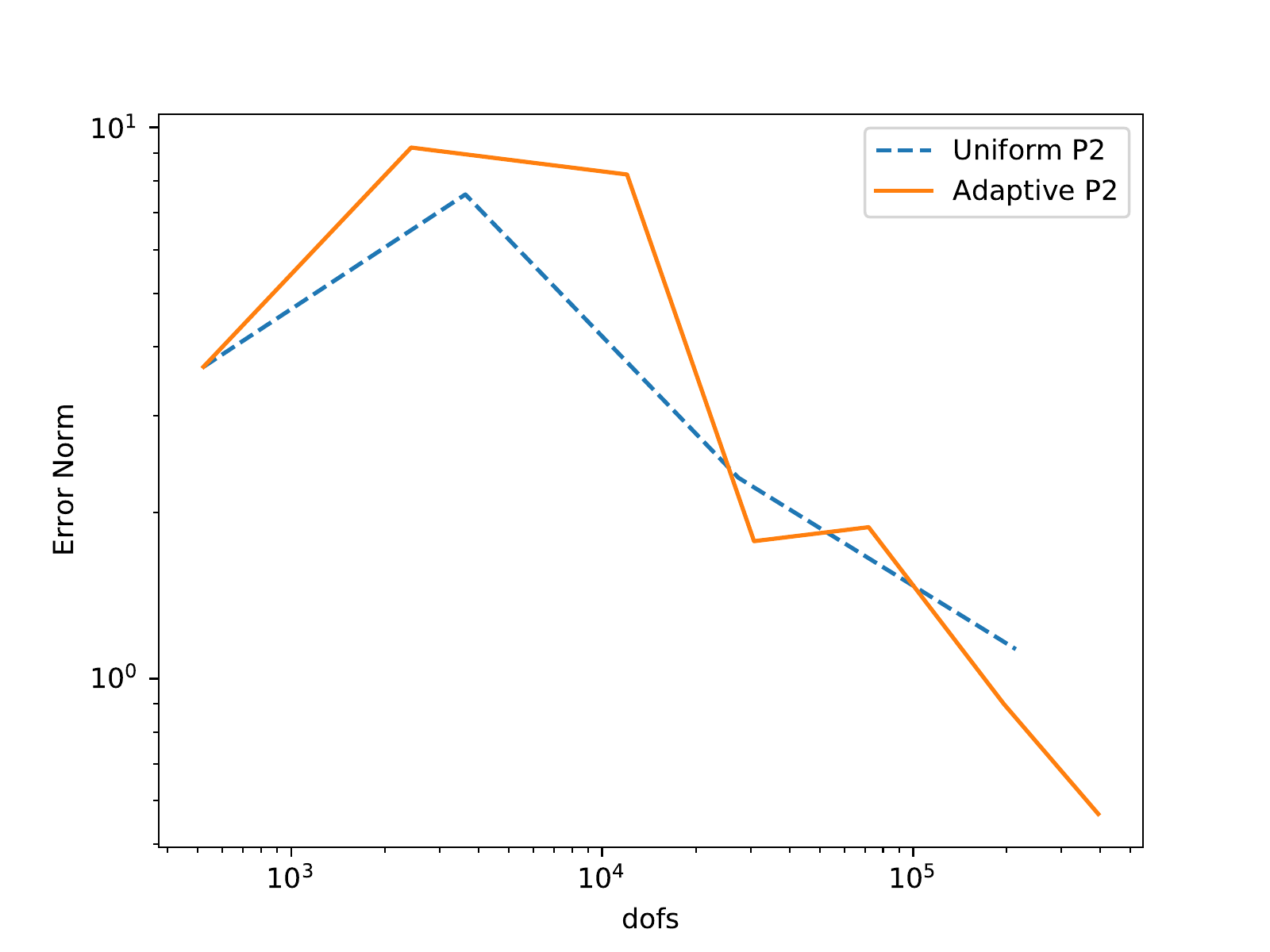}}
  \subfigure[ \label{fig:trans_comp_r_gradu_p2}  Gradient and flux approximation comparison for adaptive refinements.]{\centering 
 \includegraphics[width=0.44\textwidth]{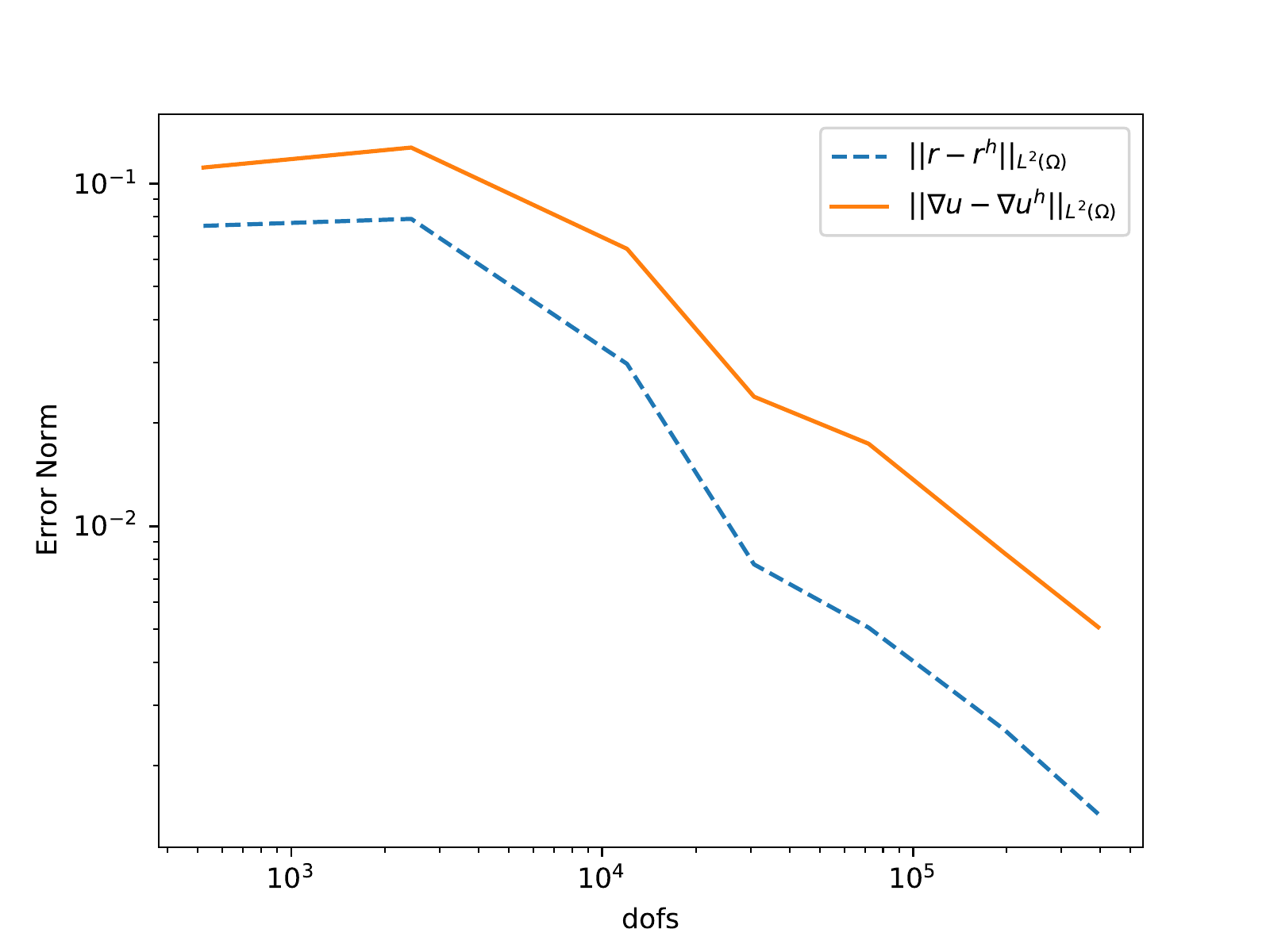}}
  \caption{\label{fig:trans_convegence_results_CH2} Convergence results for the transient Cahn-Hilliard 
  problem for quadratic basis functions.}
\end{figure}
\textcolor{black}{In Figures~\ref{fig:trans_comp_energy_p1} and~\ref{fig:trans_comp_energy_p2} the convergence 
histories in terms of the total norm on $\UUUT$ for linear and quadratic approximations are shown. The
rates of convergence of the uniform  refinements are $h^p$ and $h^{p-0.1}$, for linear and 
quadratic approximations, respectively. Since this norm contains the $H^1$ norm of $u$, we  
expect $h^p$ convergence, as indicated in~\eqref{eq:conv_rates}. 
The preasymptotic range of convergence ends at roughly $60,000$ and $80,000$ degrees of freedom for linear and quadratic approximations, respectively, at which point 
the adaptive refinements become superior. In the quadratic case, the difference between uniform and 
adaptive refinements is less pronounced in this preasymptotic range as the 
second order polynomials better approximate the sinusoidal exact solution. 
In Figures~\ref{fig:trans_comp_r_gradu_p1} and~\ref{fig:trans_comp_r_gradu_p2} we compare the errors 
in the flux variable $\rr$ and the gradient $\Nabla u$ for the case of adaptive mesh refinements. 
For both degrees of approximation we consider, the error in the flux variable is significantly lower than in the 
gradient. Compared to the stationary case in Figure~\ref{fig:comp_r_gradu}, the effect is more
pronounced for the transient problem. We attribute this to the convective nature of the time derivative 
term which has the greatest benefits of the stability property of the AVS-FE approximations. }

\subsection{Phase Transition Problem} 
\label{sec:2d_prob}
\begin{figure}[h]
{\centering
\hspace{0.1in} \includegraphics[width=0.65\textwidth]{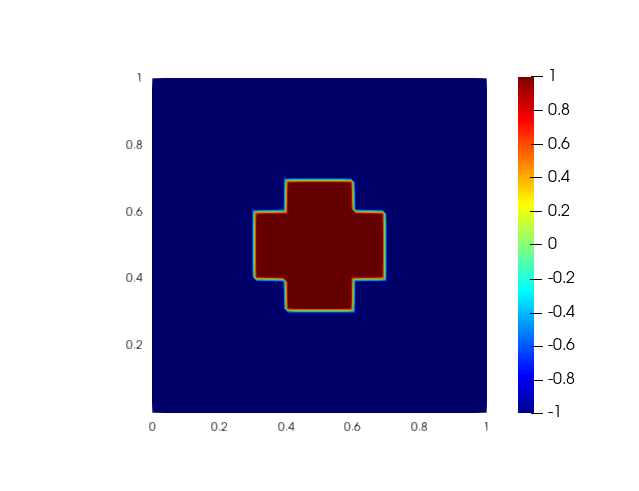}}
  \caption{\label{fig:2d_ex_ic} Initial condition for the two-dimensional model problem.}
\end{figure}

In this 
section we consider a two-dimensional benchmark problem for the Cahn-Hilliard equation 
as the target physical application of mineral separation falls into this category.
This commonly applied problem for the Cahn-Hilliard equation 
governs the evolution of two distinct phases in a medium, see, 
e.g.,~\cite{goudenege2012high,brenner2020robust}. The problem is chosen as it depicts 
a phase transformation and convergence towards a steady state.  
We consider physical properties as chosen by Brenner \emph{et al.}~\cite{brenner2020robust}: $D = 1$, $\lambda = 0.01$, and the spatial domain consist of the unit square. 
The domain is initially occupied by two phases of material, one of which is shaped like 
a cross. Inside the cross, the phase is given the value $+1$ whereas it is $-1$ outside
the cross. This initial condition is given by the piecewise constant function 
%
%
%
%
shown in Figure~\ref{fig:2d_ex_ic}. Based on 
this initial condition, the boundary conditions are $u = -1$ and $q = 0$ on $\partial \Omega$. 

\begin{figure}[h!]
\subfigure[ \label{fig:50tau} $t = 0.0015625s$. ]{\centering
 \includegraphics[width=0.45\textwidth,angle=-90,origin=c]{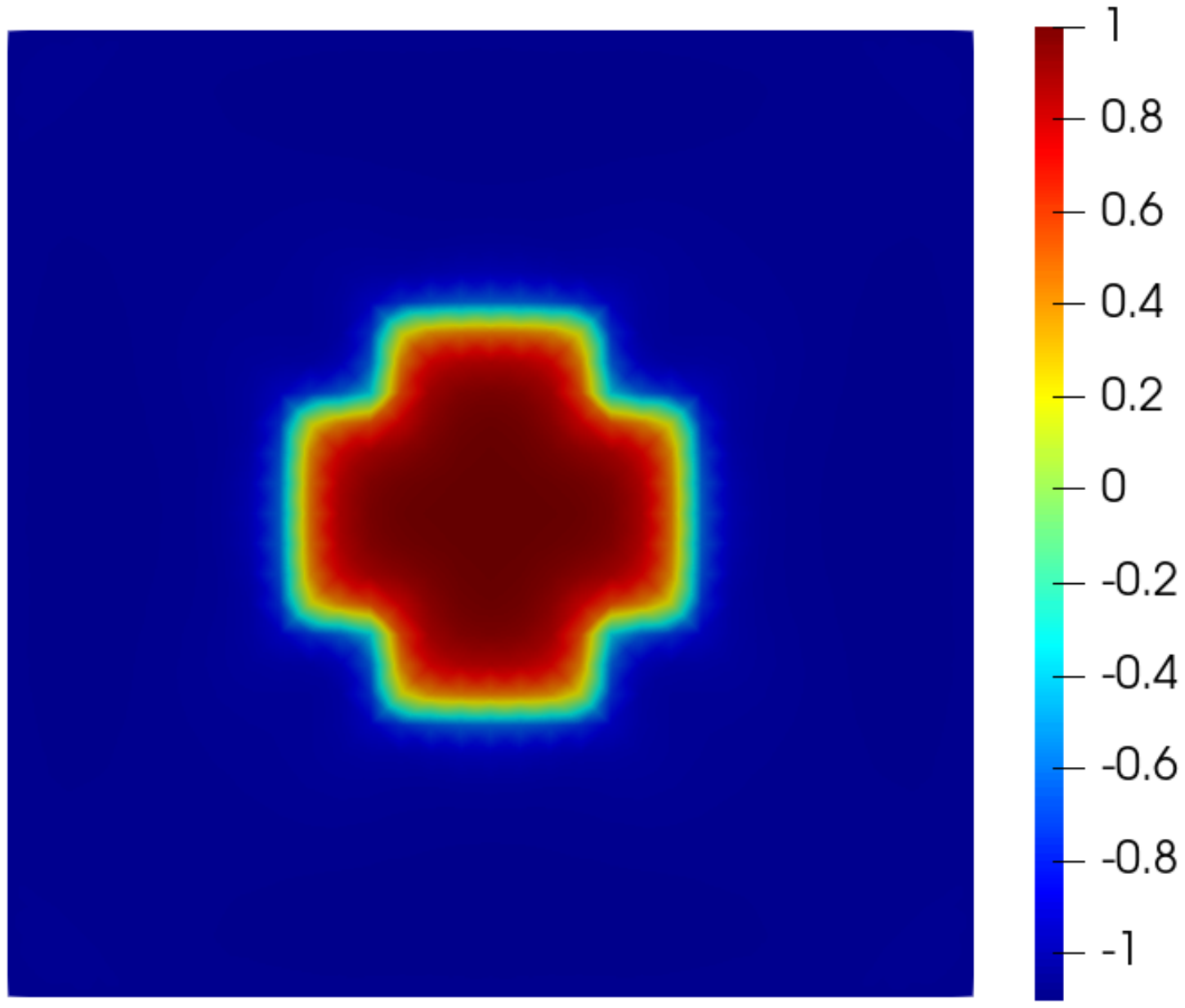}}
  \subfigure[ \label{fig:150tau} $t = 0.0046875s$.]{\centering
 \includegraphics[width=0.45\textwidth,angle=-90,origin=c]{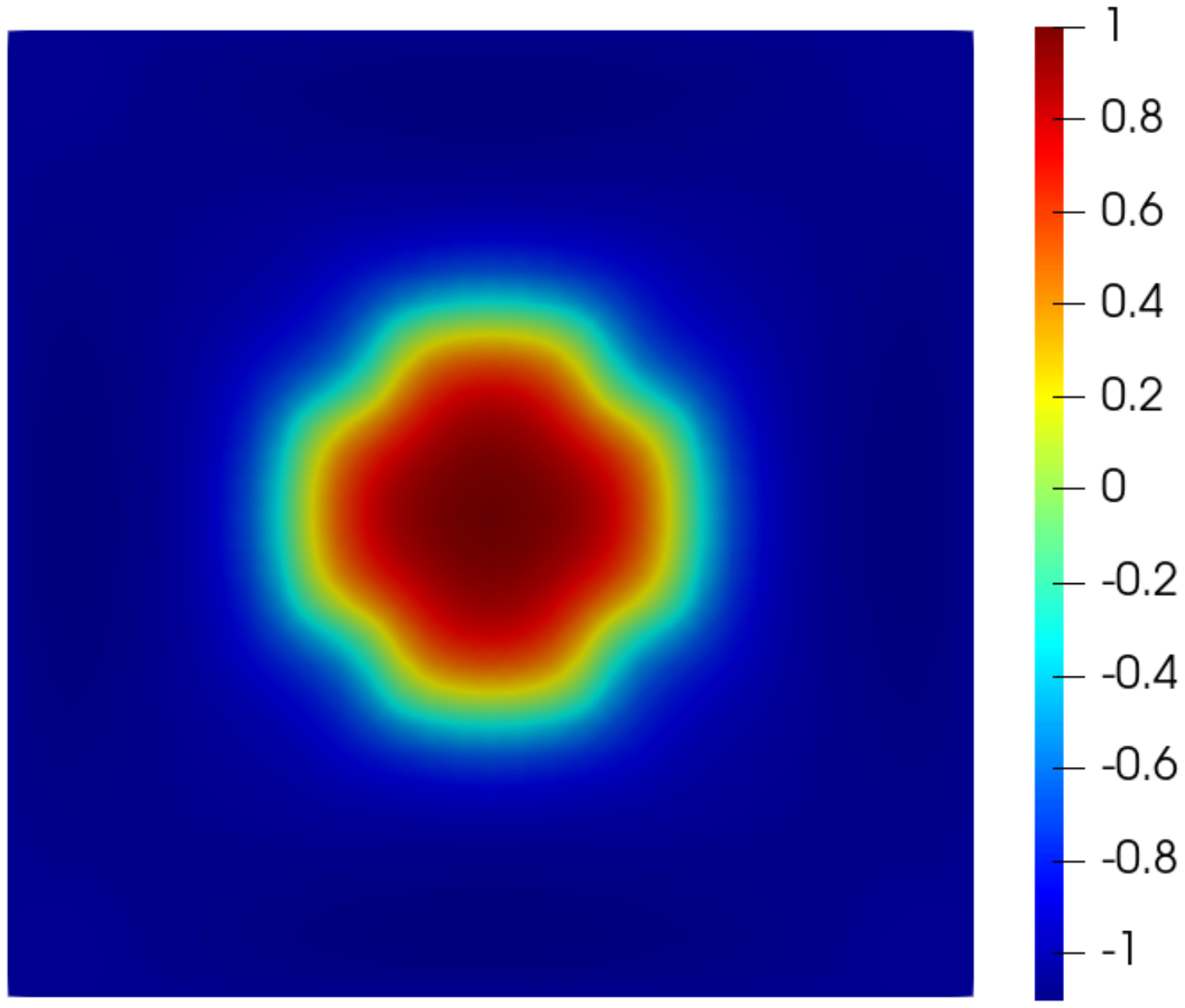}}
  \caption{\label{fig:50150tau} AVS-FE approximation $u^h$ of Cahn-Hilliard equation with initial condition shown in Figure~\ref{fig:2d_ex_ic} }
\end{figure}

\begin{figure}[h!]
\subfigure[ \label{fig:300tau}  $t = 0.09375s$. ]{\centering
 \includegraphics[width=0.45\textwidth,angle=-90,origin=c]{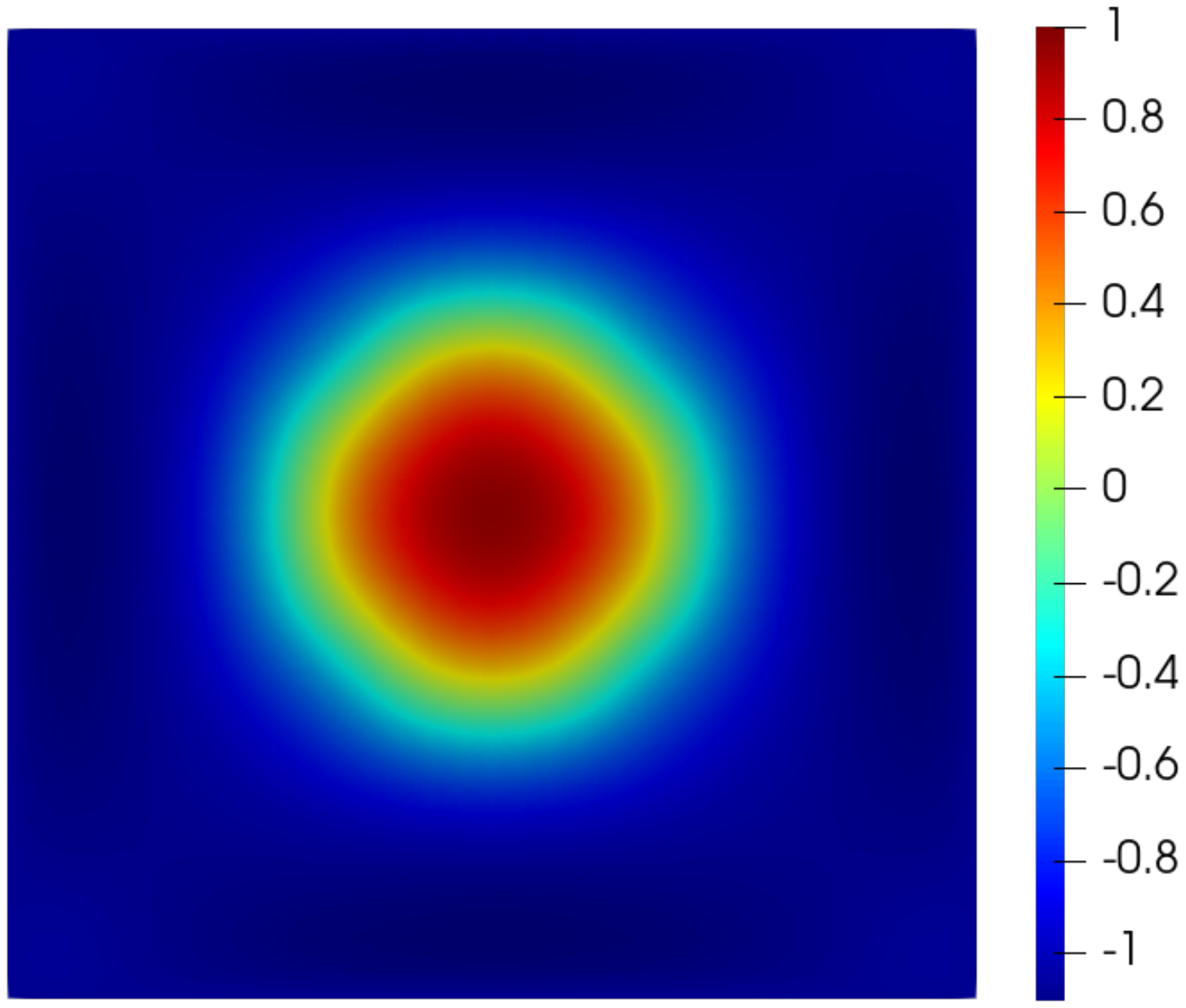}}
  \subfigure[ \label{fig:500tau} $t = 0.015625s$.]{\centering
 \includegraphics[width=0.45\textwidth,angle=-90,origin=c]{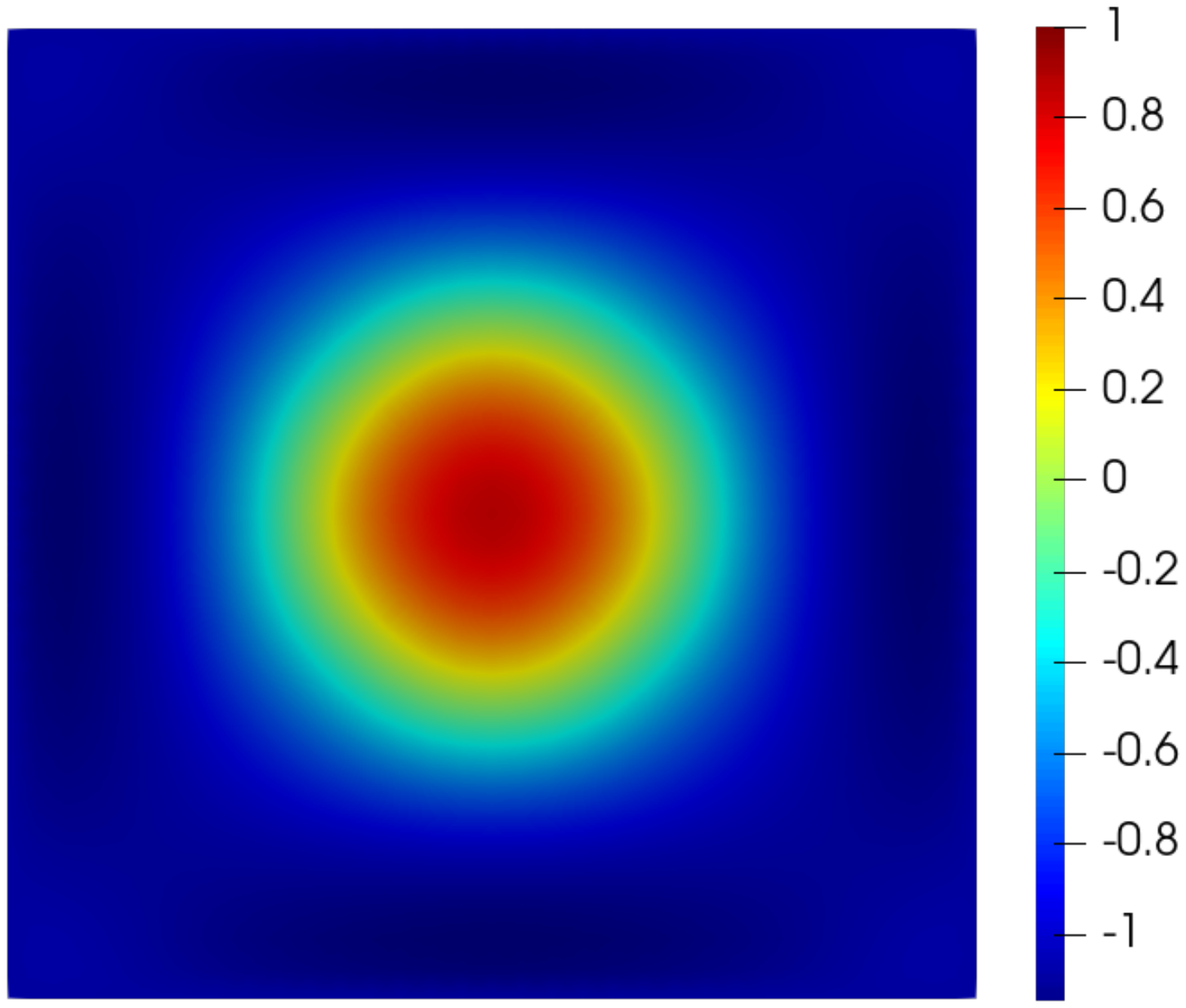}}
  \caption{\label{fig:300500tau} AVS-FE approximation $u^h$ of Cahn-Hilliard equation with initial condition shown in Figure~\ref{fig:2d_ex_ic} }
\end{figure}

To facilitate \textcolor{black}{visual} comparison with the results presented by Brenner \emph{et al.}~\cite{brenner2020robust}, we consider a case in which 
the final time $T = 0.015625s$, as the binary mixture is expected to have reached a 
steady state at this time. The computations are performed \textcolor{black}{on a fixed uniform mesh} by employing the "extruded mesh" feature 
of Firedrake~\cite{rathgeber2017firedrake}. To this end, we consider the case in which 
$\Omega_T$ is discretized by $64\times64$ triangular elements spatially that 
are extruded into triangular prism elements in the time domain of width equal to half the final time. 
We use linear polynomials spatially and fifth order polynomials in time for 
all \textcolor{black}{continuous and discontinuous} trial variables in this case. \textcolor{black}{Other choices of mesh partitions and approximation orders are possible as the method remains stable for any choice of $h$ and $p$. However, this choice is based upon extensive numerical experimentation as it provides good accuracy for a coarse mesh partition in 
time consisting of only two elements. }  
In Figures~\ref{fig:50150tau} and~\ref{fig:300500tau}, the solution is shown at $t = 0.0015625s, 0.0046875s, 0.09375s$ and $0.015625s$, 
respectively. The transformation of the binary phases from the initial to the steady state proceeds 
as expected and at the final time has reached the steady state in which the two 
phases are separated by a circle.

Using the AVS-FE method in both space and time makes the comparison to other FE methods 
for the Cahn-Hilliard equation \textcolor{black}{using the method of lines} non-trivial \textcolor{black}{and inappropriate since the computational cost is distributed in a completely different fashion.}
 In terms of accuracy, the AVS-FE method delivers 
comparable results to those reported in, e.g.,~\cite{goudenege2012high,brenner2020robust} based on visual inspection and 
comparison of the results.
 The consideration of computational 
complexity is postponed to future research efforts in which other methods
are to be considered for the time discretization.  
For this particular problem, the total number of degrees of freedom is $2,094,797$, while 
the number of degrees of freedom corresponding to the solution variables 
$(u^h,q^h,\rr^h,\ttt^h)$ is $371,800$.

\begin{figure}[h!]
     \centering
  \hspace*{-0.26in}   \scalebox{.6}{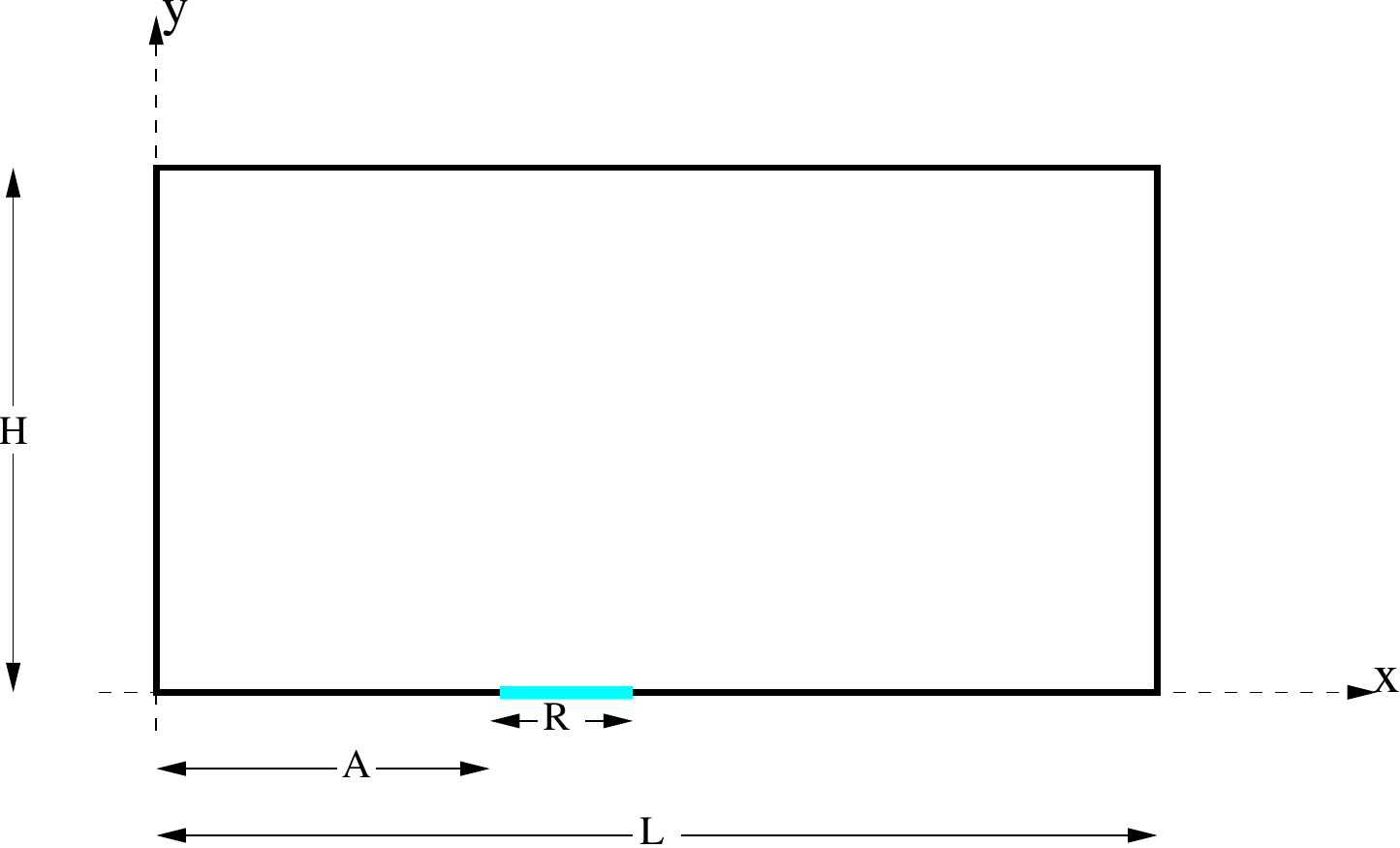 }%
  \caption{\label{fig:2D} 2D mineral separation model domain.}
\end{figure}

\textcolor{black}{
\subsection{Mineral Separation} }
\label{se:mineral_ref}
\textcolor{black}{ The development of the space-time AVS-FE method for the Cahn-Hililiard equation was 
to model the mineral separation process described in Section~\ref{sec:introduction}.   Thus, here
 we consider the application of the Cahn-Hilliard equation for the analysis of  a
mineral separation experiment. 
The experimental setup consists of a closed box into which a mineral powder is introduced 
and the separation process takes place on a substrate (see blue section in Figure~\ref{fig:2D}) 
which has 
been treated to attract the desirable mineral particles.
The accumulated mineral on the substrate disk is then collected in an appropriate fashion and for the sake 
of simplicity, we consider only the mineral accumulation in this model. 
The computational domain is a cross section of the experimental separator and is
rectangular with part of its boundary being the 
treated disk,  as shown in Figure~\ref{fig:2D}.  The physical dimensions shown in this figure are: H $= 0.3048m$, L $= 0.6604m$, A $= 0.254m$, and R $= 0.1397m$, i.e. $\Omega = (0,0.6604m)\times(0,0.3048m)$. 
Particularly, the region R in Figure~\ref{fig:2D} represents the substrate and is the location 
of the mineral accumulation.}

\textcolor{black}{
In this heuristic model we pick parameters D = 1 and $\lambda$ = 0.1. Note that the  proper physical 
parameters are to be estimated using an inverse finite element process using experimental data  
from the experiment in the ongoing design process.
 The initial condition is a concentration of 0 throughout the domain $\Omega$ and to model the 
 buildup of minerals, we employ 
 a Dirichlet boundary condition with a mineral concentration equal to 1 over the disk region. 
The remainder of the boundary $\partial \Omega$ is considered to have zero mineral concentration.
We implement this problem in the same fashion as the preceding verification
in Firedrake~\cite{rathgeber2017firedrake} and use} \textcolor{black}{a fixed uniform mesh. } \textcolor{black}{
Hence, the mesh partition consists of $64\times64$ triangular elements extruded into two space-time 
triangular prisms and the basis functions are 
polynomials that are linear in space and fourth order in time.}

\begin{figure}
     \centering
     \includegraphics[trim={5cm 0cm 0cm 9cm},clip,width=.7\textwidth]{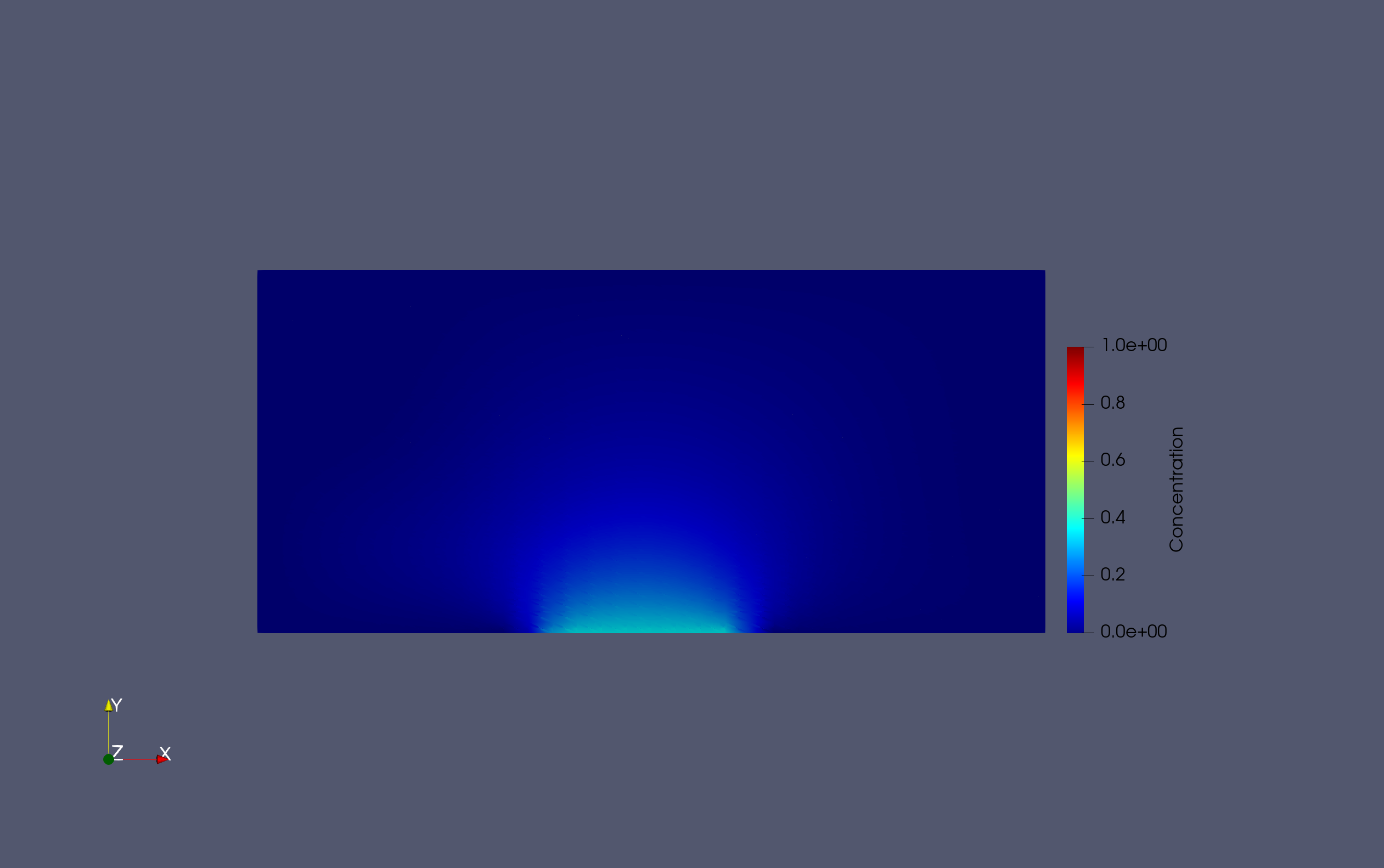}
  \caption{\label{fig:InitialCond} Mineral separation process at $t = 0.001s$.}
\end{figure}

\begin{figure}
     \centering
     \includegraphics[trim={5cm 0cm 0cm 9cm},clip,width=.7\textwidth]{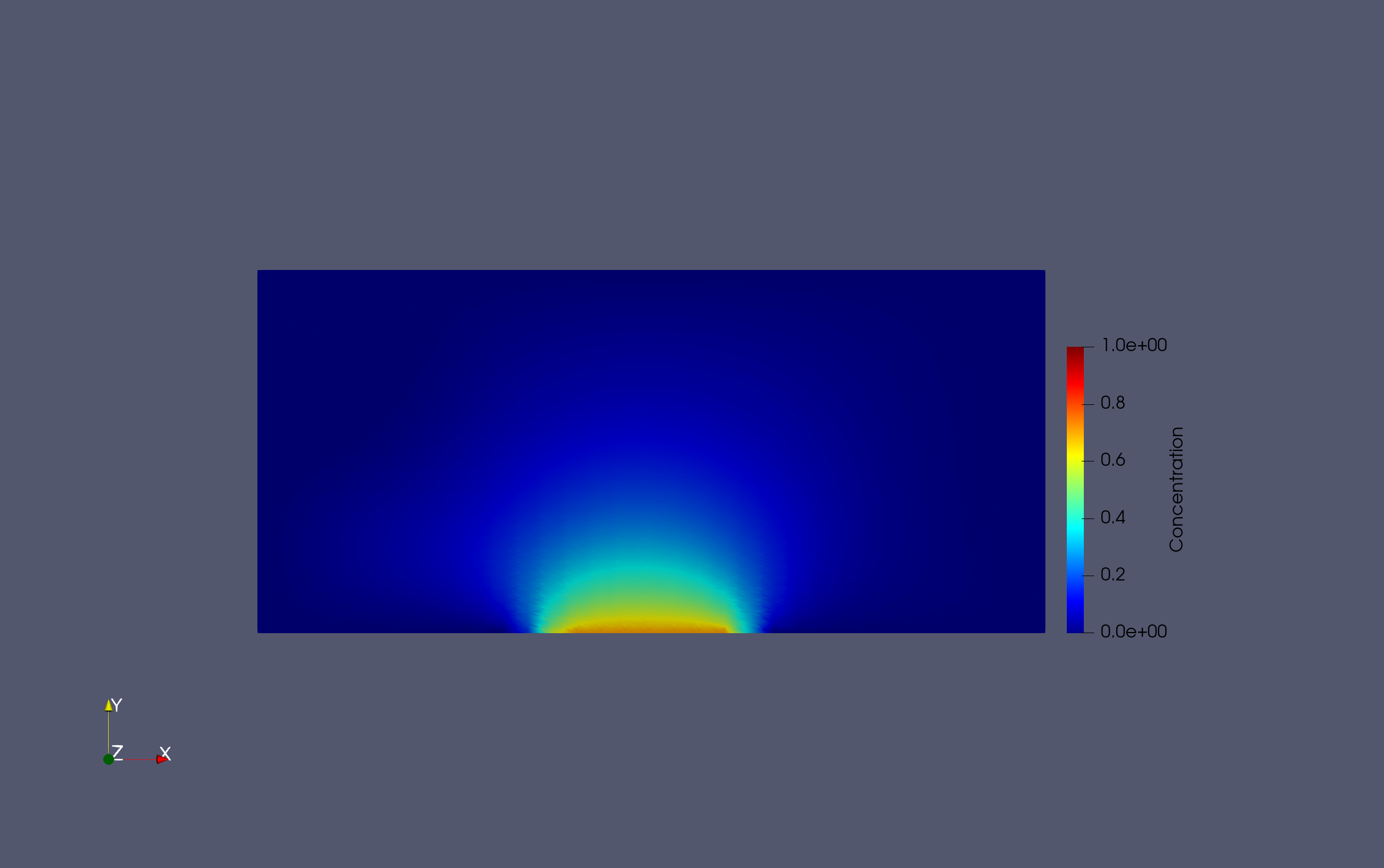}
  \caption{\label{fig:InitialCond2} Mineral separation process at $t = 0.002s.$}
\end{figure}

\begin{figure}
     \centering
      \includegraphics[trim={5cm 0cm 0cm 9cm},clip,width=.7\textwidth]{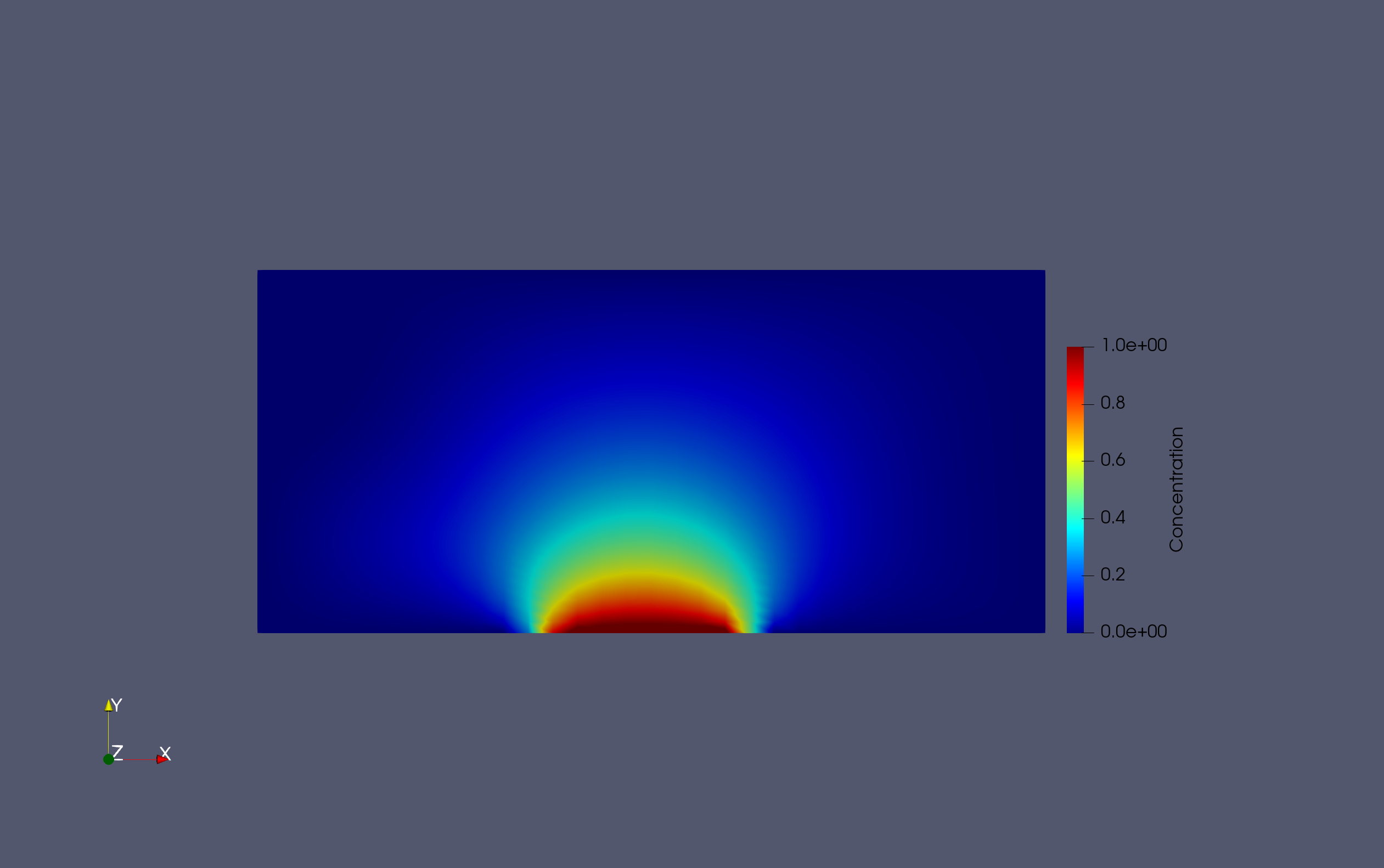}
  \caption{\label{fig:InitialCond3} Mineral separation process at $t = 0.004s.$}
\end{figure}

\textcolor{black}{ In
Figures~\ref{fig:InitialCond},~\ref{fig:InitialCond2}, and~\ref{fig:InitialCond3} the buildup of mineral
is shown at times $t = 0.001s$, $0.002s$, and $0.004s$ respectively. 
As expected, the mineral layer grows vertically since we do not incorporate effects of airflow in this case. 
Hence, we conclude that the Cahn-Hilliard equation is an appropriate model for the buildup of material in the proposed mineral processing application.}

\section{Conclusions}
\label{sec:conclusions}
We have presented an application of the AVS-FE method to the Cahn-Hilliard equation 
to establish stable AVS-FE approximations of the 
Cahn-Hilliard equation in both space and time.

The AVS-FE method results in FE approximations that converge to the exact solution at optimal rates. This was illustrated by \textcolor{black}{both stationary and transient two-dimensional verifications to which manufactured solutions exist in Section~\ref{sec:1d_prob} }. 
Furthermore, the AVS-FE method comes with a built-in error estimate and the resulting error indicator has been used to successfully drive mesh adaptive refinements in both space and time.  
\textcolor{black}{As expected, the adaptive mesh refinements deliver lower FE approximation errors in 
appropriate norms than uniform refinements.}
\textcolor{black}{A benchmark} problem from literature was also implemented in which the Cahn-Hilliard 
problem leads to a phase transition between two binary phases in a medium. 
The results for this case show that the AVS-FE method is capable of delivering 
accurate space-time computations for spatially two-dimensional problems \textcolor{black}{requiring only 
a single nonlinear global solve}. 
\textcolor{black}{Finally, in an effort to verify and Cahn-Hilliard equations  as a model for the 
mineral separation application of interest we consider a heuristic model problem governing the 
growth of a mineral in a simplified mineral separator. }

Focus is given here to the use of the AVS-FE method in both space and time due to its 
discrete stability. Other time discretizations can also be considered, e.g., 
finite difference techniques or generalized $\alpha$ 
method~\cite{deng2019high,chung1993time}. 
For now we leave these techniques for future works that involve spatially 
three-dimensional problems, in which the space time approach of this paper becomes less feasible due to the need for four-dimensional mesh generation. 
While the computational cost of the space-time AVS-FE approach we have presented in 
this paper is high \textcolor{black}{compared to existing techniques} it has the  
 advantage that we it solve transient BVPs using 
a single global solve in space and time. Additionally, the stability property of the AVS-FE 
enables users to start with a very coarse space-time mesh which can be adaptively refined 
using the built-in error indicator. Hence, our 
method is a viable option to existing time stepping algorithms, and the 
 entire solution process is easily amendable to parallel processing. 
Furthermore, a static condensation process through a Schur complement can further 
reduce the computational cost. While we have not implemented this in our current computational framework it is to be one of our foci for future research 
efforts \textcolor{black}{as well as implementation into other DPG FE software such as Camellia \cite{roberts2019camellia}}. \textcolor{black}{The existence of a Fortin operator~\cite{nagaraj2017construction,demkowicz2020construction} is at this point based on conjecture from numerical experiments and future research will attempt a mathematical construction of these operators for the AVS-FE method.}
\textcolor{black}{Currently, the development of an inverse FE process to estimate the physical parameters $D$ and $\lambda$ for the mineral separation application is ongoing.}

\section*{Acknowledgements}
The authors are grateful for the contributions of Professor Jon Kellar,  Professor William Cross, and 
Mr. Bernardo Sansao of the Department of Materials and Metallurgical Engineering at the SDSM\&T through fruitful discussions on mineral separation and the Cahn-Hilliard equation.

This work has been supported by the NSF CBET Program,
under  NSF Grant Number 1805550.

\bibliographystyle{elsarticle-num}
 \bibliography{references_eirik}

 \newpage
\appendix 
\renewcommand{\theequation}{A.\arabic{equation}}
\renewcommand{\thesection}{A.\arabic{section}}

\color{black}{
\section*{Appendix A}
\label{app:append_A}
To establish the well-posedness of the AVS-FE weak formulation~\eqref{eq:weak_form_CH} we apply 
the Babu{\v{s}}ka Lax-Milgram Theorem~\cite{babuvska197finite}. To accomplish this, we follow the 
steps of Carstensen, Demkowicz and Gopalakrishnan~\cite{carstensen2016breaking} to show 
that the stability of the AVS-FE weak formulation is inherited from its unbroken counterpart.
Hence, we proceed by first showing the well-posedness of a mixed weak form of a linear 
Cahn-Hilliard problem, and subsequently showing that the corresponding weak form with a broken 
test space is well posed as a consequence. In the following, we denote by $C$ a generic mesh
independent constant.

For this analysis, we consider a 
simplified linear and stationary form of the Cahn-Hilliard BVP:
\begin{equation}  \label{eq:Cahn-Hilliard_BVP_simple}
\boxed{
\begin{array}{l}
\text{Find }  u  \text{ such that:}    
\\[0.05in] 
\qquad 
\begin{array}{rcl}
 \ds  D \, \Delta \left[ - u  
 -\lambda \, \Delta \, u \right]  & = & f, \quad \text{ in } \, \Omega, 
 \\[0.05in]
 \qquad \ds u = - u -\lambda \, \Delta \, u &  = & 0, \quad \text{ on } \,  \partial \Omega,
 \end{array}
 \end{array}
}
\end{equation}
where we have assumed that non-homogeneous Dirichlet boundary conditions are replaced by an
 appropriate source term $f \in \SLTO$. The corresponding first-order system is:
\begin{equation} \label{eq:CH_IBVP_first_order_simple}
\boxed{
\begin{array}{l}
\text{Find }  (u,q,\rr,\ttt) \text{ such that:}    
\\[0.05in] 
\qquad 
\begin{array}{rcl}
\ds   \Nabla u - \rr & =  & \bfm{0}, \quad \text{ in } \, \Omega, 
\\
\ds   - u - \lambda \, \Nabla \cdot \rr - q & =  & 0, \quad \text{ in } \, \Omega, 
  \\
  \ds   \Nabla q - \ttt & =  & \bfm{0}, \quad \text{ in } \, \Omega, 
  \\
  \ds    - D \, \Nabla \cdot \ttt & =  & -f, \quad \text{ in } \, \Omega, 
 \\[0.025in]
 \qquad u &  = & 0, \quad \text{ on } \, \partial \Omega, 
 \\
 \qquad q & = & 0, \quad \text{ on } \,  \partial \Omega. 
 \end{array}
 \end{array}
}
\end{equation}
Multiplying this first-order system by test functions $(v,\sss,w,\pp)$ and applying integration by parts
to the terms $\Nabla \cdot \rr$ and $\Nabla \cdot \ttt$ and a strong enforcement of boundary conditions 
yields the weak formulation:
\begin{equation} \label{eq:weak_form_CH_app}
\boxed{
\begin{array}{ll}
\text{Find } (u,q,\rr,\ttt) \in \UUUA  \text{ such that:}
\\[0.05in]
   \quad b((u,q,\rr,\ttt),(v,w, \sss, \pp)) = f(w), \quad \forall (v,w, \sss, \pp)\in \VOM,
 \end{array}}
\end{equation}
where $\UUUA = \VOM = H^1_0(\Omega)  \times H^1_0(\Omega) \times \SLTOvec \times \SLTOvec$ and:
\begin{equation} \label{eq:B_and_F_CH_app}
\begin{array}{l}
b((u,q,\rr,\ttt),(v,w, \sss, \pp)) \isdef 
\ds  \int_{\Omega}\left[  \, ( \Nabla u - \rr) \cdot \sss \, + ( \Nabla q - \ttt ) \cdot \pp \, - (  u + q  ) \, v 
+ \ds  \, \lambda \,   \rr \cdot \Nabla v  +  D \, \ttt \cdot \Nabla  w \right] \dx,
 \\[0.1in]
 \hspace{1in} \ds f(w) \isdef \int_{\Omega} -f \, w \dx.
 \end{array}
\end{equation}
By defining  two bilinear forms:
\begin{equation} \label{eq:a_and_c_CH_app}
\begin{array}{l}
a((u,q,\rr,\ttt),(v,w, \sss, \pp)) \isdef 
\ds  - \int_{\Omega}\left[  \, uv + qv + \rr \cdot \sss  + \ttt \cdot \pp \right] \dx,
 \\[0.1in]
 \hspace{0.4in} \ds c((v,w),(\eps, \boldphi)) \isdef \int_{\Omega} \eps \cdot \Nabla v + \boldphi \cdot \Nabla w \dx,
 \end{array}
\end{equation}
we can recast~\eqref{eq:B_and_F_CH_app} as:
\begin{equation} \label{eq:mixed_CH_app}
\boxed{
\begin{array}{l}
\text{Find } (u,q,\rr,\ttt) \in \UUUA  \text{ such that:}
\\[0.05in]
a((u,q,\rr,\ttt),(v,w, \sss, \pp)) +  c((v,w),(\lambda \rr, D\, \ttt)) = f(w),
\\[0.1in]
\hspace{2.9in} \quad \forall (v,w, \sss, \pp)\in \VOM.
 \\[0.1in]
  c((u,q),( \sss, \pp)) \hspace{1.6in}= 0, 
 \end{array} }
\end{equation}
\begin{prp} 
\label{prp:brezzi_cond}
Let $v \in H^1_0(\Omega)$.  Then, $c_1(v,\eps) = \int_{\Omega} \eps \cdot \Nabla v  \dx$ satisfies the inf-sup condition:
\begin{equation} \label{eq:small_inf}
\begin{array}{l}
\exists \gamma >0 : \supp{v \in H^1(\Omega) } \ds 
     \frac{|c_1(v,\eps)|}{\norm{v}{\SHOO}} \ge \gamma \norm{\eps}{\SLTO}
 \end{array}
\end{equation}
\end{prp} 
\emph{Proof}: Because $v \in H^1_0(\Omega)$ the Poincar\' e inequality gives:
\begin{equation*} 
\begin{array}{l}
  \supp{v \in H^1(\Omega) } \ds \frac{| \int_{\Omega} \eps \cdot \Nabla v  \dx |}{\norm{v}{\SHOO}} \ge 
  \supp{v \in H^1(\Omega) } \ds \frac{| \int_{\Omega} \eps \cdot \Nabla v  \dx |}{ C \, \norm{\Nabla v}{\SLTO}}. 
 \end{array}
\end{equation*}
Next, we pick $\Nabla v = \eps$ to get:
\begin{equation*} 
\begin{array}{l}
  \supp{v \in H^1(\Omega) } \ds \frac{| \int_{\Omega} \eps \cdot \Nabla v  \dx |}{\norm{v}{\SHOO}} \ge 
   \ds \frac{| \int_{\Omega} \eps \cdot \eps \dx |}{ C \, \norm{\eps}{\SLTO}} =   C \, \norm{\eps}{\SLTO}, 
 \end{array}
\end{equation*}
and the proof is complete with $\gamma = C$.
\newline \noindent ~\qed

\begin{lem}
\label{lem:well_posed_cont_APP}
Let $(u,q,\rr,\ttt) \in \UUU$. Then, the mixed weak problem~\eqref{eq:mixed_CH_app} is well posed.
\end{lem} 
\emph{Proof}: Since the form $a(\cdot,\cdot)$ is clearly coercive, then
the mixed problem~\eqref{eq:mixed_CH_app} is well posed because all terms of $c(\cdot,\cdot)$ can be 
shown to satisfy an  
inf-sup condition shown in Proposition~\ref{prp:brezzi_cond}, i.e., Brezzi's condition 
(See Theorem II.1.1 in~\cite{BrezziMixed}).
\newline \noindent ~\qed

\noindent Lemma~\ref{lem:well_posed_cont_APP} leads to an inf-sup condition (see, e.g.,~\cite{demkowicz2020construction}) of the form:
\begin{equation} \label{eq:inf-supp}
\begin{array}{l}
\exists \gamma >0 : \supp{(v,w, \sss, \pp) \in \VOM } \ds 
     \frac{|b((u,q,\rr,\ttt),(v,w, \sss, \pp))|}{\norm{(v,w, \sss, \pp)}{\VOM}} \ge \gamma \, \norm{(u,q,\rr,\ttt)}{\UUUA}.
 \end{array}
\end{equation}

Before showing the well-posedness of the AVS-FE weak formulation, let us use the continuous 
bilinear form $b(\cdot,\cdot)$ to write the corresponding AVS-FE weak form:
\begin{equation} \label{eq:AVS-weak_app}
\begin{array}{l}
B((u,q,\rr,\ttt),(v,w, \sss, \pp)) = b((u,q,\rr,\ttt),(v,w, \sss, \pp)) + \dualp{D\,\ttt\cdot\nn}{w_m}_{\Gamma_h}+\dualp{\lambda\,\rr\cdot\nn}{v_m}_{\Gamma_h},
 \end{array}
\end{equation}
where $\dualp{D\,\ttt\cdot\nn}{w_m}_{\Gamma_h}+\dualp{\lambda\,\rr\cdot\nn}{v_m}_{\Gamma_h} \isdef \summa{\Kep}{}
\oint_{\dKm} \{ D \,  \gamma^m_\nn(\ttt) \, \gamma^m_0(w_m) + \lambda \,  \gamma^m_\nn(\rr) \, 
\gamma^m_0(v_m) \, \} \dss$. 
We also require the following intermediate results.
\begin{prp} 
\label{prp:skel_infsup}
Let $\ttt,\rr \in \SHdivP$ and $v_m,w_m \in \SHOP$.  Then:
\begin{equation} \label{eq:inf_sup_skel}
\begin{array}{l}
\exists \gamma^s >0 : \supp{(v,w) \in \SHOP\times\SHOP } \ds 
     \frac{|\dualp{D\,\ttt\cdot\nn}{w_m}_{\Gamma_h}+\dualp{\lambda\,\rr\cdot\nn}{v_m}_{\Gamma_h}|}{\norm{(v,w)}{\SHOP}} \ge \gamma^s \, \norm{(\rr,\ttt)}{\UUUhat},
 \end{array}
\end{equation}
where $\SHdivP$ denotes the broken $H(div)$ space and $\norm{(\rr,\ttt)}{\UUUhat}$ is the minimum 
energy extension norm:
\begin{equation} \label{eq:min_en_ext}
\begin{array}{l}
\norm{(\rr,\ttt)}{\UUUhat}  \isdef \summa{\Kep}{}
\oint_{\dKm} \{   \gamma^m_\nn(\ttt) \, \gamma^m_0(w_m) + \gamma^m_\nn(\rr) \, 
\gamma^m_0(v_m) \, \} \dss= \,\rm{inf} \,( \norm{\rr}{\SHdivO}^2 + \norm{\ttt}{\SHdivO}^2 )^{1/2}
 \end{array}
\end{equation}
\end{prp} 
\emph{Proof}: see Theorem 2.3 in~\cite{carstensen2016breaking}.
\newline \noindent ~\qed
$
$
\begin{prp} 
\label{prp:zero_skel_trac}
Let $\ttt,\rr \in \SHdivO$ and $v_m,w_m \in \SHOP$.  Then:
\begin{equation} \label{eq:zero_traces}
\begin{array}{l}
\dualp{D\,\ttt\cdot\nn}{w_m}_{\Gamma_h}+\dualp{\lambda\,\rr\cdot\nn}{v_m}_{\Gamma_h} = 0, \quad \forall v_m,w_m \in \SHOO.
 \end{array}
\end{equation}
\end{prp} 
\emph{Proof}: Since functions in $\SHdivO$ have zero jump across mesh interfaces both terms must vanish for all single valued functions $v_m$ and $w_m$ on $\Gamma_h$, see Theorem 2.3 in~\cite{carstensen2016breaking}.
\newline \noindent ~\qed

\noindent Now, Lemma~\ref{lem:well_posed_cont_APP} and Propositions~\ref{prp:skel_infsup} 
and~\ref{prp:zero_skel_trac}
correspond to the necessary assumptions of Theorem 3.3 in~\cite{carstensen2016breaking}. Hence, 
we replicate their arguments  to show the following assertion:
\begin{lem}
\label{lem:well_posed_broken_APP}
Let $(u,q,\rr,\ttt) \in \UUU$. Then, the AVS-FE weak formulation~\eqref{eq:AVS-weak_app}
is well posed.
\end{lem} 
\emph{Proof}: Since continuity of the bilinear form and linear functional can be shown in a 
straightforward manner by the Cahchy-Schwarz inequality and successive integration by parts. 
The inf-sup condition is the last required point of the 
Babu{\v{s}}ka Lax-Milgram Theorem~\cite{babuvska197finite} for the well-posedness of the weak
formulation. By~\eqref{eq:inf-supp} we have:
\begin{equation*}
\begin{array}{l}
\ds C \, \norm{(u,q,\rr,\ttt)}{\UUUA} \le  \supp{(v,w, \sss, \pp) \in \VOM } \ds 
     \frac{|b((u,q,\rr,\ttt),(v,w, \sss, \pp))|}{\norm{(v,w, \sss, \pp)}{V}},
 \end{array}
\end{equation*}
notice that in the denominator we have replaced $\norm{(v,w, \sss, \pp)}{\VOM}$ with
$\norm{(v,w, \sss, \pp)}{V}$, which are equivalent for functions $(v,w, \sss, \pp) \in \VOM$.
We then add the duality pairings, i.e, zero (see Proposition~\ref{prp:zero_skel_trac}) in the numerator:
\begin{equation*}
\begin{array}{l}
\ds C \, \norm{(u,q,\rr,\ttt)}{\UUUA} \le  \supp{(v,w, \sss, \pp) \in \VOM } \ds 
     \frac{|b((u,q,\rr,\ttt),(v,w, \sss, \pp))+\dualp{D\,\ttt\cdot\nn}{w_m}_{\Gamma_h}+\dualp{\lambda\,\rr\cdot\nn}{v_m}_{\Gamma_h}|}{\norm{(v,w, \sss, \pp)}{V}},
 \end{array}
\end{equation*}
finally, we test with a larger space $\VV \supset \VOM$ to establish the inf-sup condition:
\begin{equation} \label{eq:int_inf_supp}
\begin{array}{l}
\ds C \, \norm{(u,q,\rr,\ttt)}{\UUUA} \le  \supp{(v,w, \sss, \pp) \in \VV } \ds 
     \frac{|B((u,q,\rr,\ttt),(v,w, \sss, \pp))|}{\norm{(v,w, \sss, \pp)}{\VV}}.
 \end{array}
\end{equation}
Next, by Proposition~\ref{prp:skel_infsup}:
\begin{equation*} 
\begin{array}{l}
\ds 
 \gamma^s \, \norm{(\rr,\ttt)}{\UUUhat} \le  \supp{(v,w, \sss, \pp) \in \VV }   \frac{|\dualp{D\,\ttt\cdot\nn}{w_m}_{\Gamma_h}+\dualp{\lambda\,\rr\cdot\nn}{v_m}_{\Gamma_h}|}{\norm{(v,w, \sss, \pp)}{V}}.
 \end{array}
\end{equation*}
By~\eqref{eq:AVS-weak_app}, we get:
\begin{equation*} 
\begin{array}{l}
\ds 
 \gamma^s \, \norm{(\rr,\ttt)}{\UUUhat} \le  \supp{(v,w, \sss, \pp) \in \VV }   \frac{|B((u,q,\rr,\ttt),(v,w, \sss, \pp)) - b((u,q,\rr,\ttt),(v,w, \sss, \pp))|}{\norm{(v,w, \sss, \pp)}{V}},
 \end{array}
\end{equation*}
which can be further bound using~\eqref{eq:inf-supp}:
\begin{equation*} 
\begin{array}{l}
\ds 
 \gamma^s \, \norm{(\rr,\ttt)}{\UUUhat} \le  \gamma \, \norm{(u,q,\rr,\ttt)}{\UUUA} + \supp{(v,w, \sss, \pp) \in \VV }    \frac{|B((u,q,\rr,\ttt),(v,w, \sss, \pp))|}{\norm{(v,w, \sss, \pp)}{V}},
 \end{array}
\end{equation*}
where the first term in the RHS can be bound by~\eqref{eq:int_inf_supp}. Finally,
we note that the norm on $\UUU$ can be expressed as:
\begin{equation*} 
\begin{array}{l}
\ds  \norm{(u,q,\rr,\ttt)}{\UUU}^2 = \norm{(\rr,\ttt)}{\UUUhat}^2 + \norm{(u,q,\rr,\ttt)}{\UUUA}^2,
 \end{array}
\end{equation*}
leading to the desired \emph{inf-sup} condition: 
\begin{equation} \label{eq:APP_inf_sup_AVS}
\begin{array}{l}
\ds C \, \norm{(u,q,\rr,\ttt)}{\UUU} \le  \supp{(v,w, \sss, \pp) \in \VV } \ds 
     \frac{|B((u,q,\rr,\ttt),(v,w, \sss, \pp))|}{\norm{(v,w, \sss, \pp)}{\VV}}.
 \end{array}
\end{equation}
\noindent ~\qed
}

\end{document}